\documentclass[12pt]{amsart}
\usepackage{amssymb,latexsym,amsmath,amsfonts}
\theoremstyle{plain}
\newtheorem{theorem}{Theorem}
\newtheorem{corollary}{Corollary}

\newtheorem{lemma}{Lemma}
\newtheorem{proposition}{Proposition}

\theoremstyle{definition}

\theoremstyle{remark}

\newtheorem{remark}{Remark}

\newcommand{\eqdef}{\stackrel{\mathrm{def}}{=}}

\def\d{\delta}
\def\a{\alpha}

\def\g{\gamma}

\def\e{\epsilon}
\def\l{\lambda}
\def\m{\mu}

\def\s{\sigma}

\def\C{\mathbb C}
\def\N{\mathbb N}
\def\R{\mathbb R}
\def\Z{\mathbb Z}

\def\lim{\operatorname{lim}}

\def\sup{\operatorname{sup}}
\def\max{\operatorname{max}}

\def\ovS{\overline{S}}

\begin{document}
\title[Von Mises statistics of a transformation]{Limit theorems for von Mises
statistics of a measure preserving transformation}
 \date{July 22, 2013}
 \author[Manfred Denker and Mikhail Gordin]{Manfred Denker
 and Mikhail Gordin }
\address {{M. Denker}\\ Mathematics Department\\ Pennsylvania State
   University\\ State College PA 16802, USA}
\email{denker@math.psu.edu}
\address{{M. Gordin}\\  V.A. Steklov Institute of Mathematics, St. Petersburg Division,
27 Fontanka emb. \\ Saint Petersburg 191023, RUSSIA\newline and
\newline
Saint Petersburg State University, Faculty of Mathematics and
Mechanics, 28 Universitetski av., Petrodvorets, Saint Petersburg
198504, RUSSIA} \email{gordin@pdmi.ras.ru}
\begin{abstract}

 For a measure preserving transformation $T$ of
a probability space $(X,\mathcal F,\mu)$ and some $d \ge 1$ we
investigate almost sure and distributional convergence of
  random variables of the form
  $$x \to \frac{1}{C_n} \sum_{0\le i_1,\dots,\,i_d<n}
  f(T^{i_1}x,\dots,T^{i_d}x),\ n=1,2, \dots,
  $$
  where  $C_1, C_2,\dots$
  are normalizing constants and the kernel $f$  belongs to an appropriate
  subspace in some $L_p(X^d\!,\, \mathcal F^{\otimes d}\!,\,\mu^d)$.  We establish
  a form of the individual ergodic theorem
for such sequences. Using a filtration compatible with $T$ and the
martingale approximation, we prove a central limit theorem in the
non-degenerate case;
for a class of canonical
(totally
degenerate) kernels and $d=2$, we also show that
  the convergence holds in distribution
  towards a quadratic form $\sum_{m=1}^{\infty}
  \lambda_m\eta^2_m$ in independent standard Gaussian variables
  $\eta_1, \eta_2, \dots$.
\end{abstract}
\maketitle

\section{Introduction} \label{1}
\subsection{Objectives and contents}
The present paper
aims to extend the theory of von Mises statistics for independent,
identically distributed random variables to the realm of strictly
stationary processes.
Every stationary process will be investigated together with a
respective measure preserving transformation of the main
probability space. Such a transformation is the only structure
used in the present article to establish  a Strong Law of Large
Numbers (SLLN) for von Mises statistics.
The Central Limit Theorem (CLT) and other weak convergence results
are treated in the framework of a filtration compatible with the
transformation. A stationary processes generating  such a
filtration will appear only in applications. It turns out that a
considerable part of the limit theory can be developed on this
basis. One of the objectives of the paper is to show that such a
relatively modest additional structure creates a suitable setting
to apply some form of the martingale approximation; indeed, the
latter is our main tool when proving the CLT-type results.
 Below, we will explain  another objective
 of the present work and its results; the latter are collected in four
statements.\\

Let $T$ be a measure preserving transformation of a probability
space $(X,\mathcal F,\mu)$. For every $d \ge 1$ and every suitable
(see the next paragraph for the elaboration) measurable function
$f: X^d \to \R$, called a \emph{kernel}, we investigate, after
normalizing appropriately, the asymptotic behavior of random
variables
\begin{equation} \label{sum}
x \mapsto  \sum_{0 \le i_1<\,n,...,\,0 \le \,i_d\,<n}
f(T^{i_1}x,...,T^{i_d}x),\, n=1,2, \dots,
  \end{equation}
  as $n$ tends to $\infty$. Every function of the form \eqref{sum}, normalized by some
  constant or not, will be called a \emph{von Mises statistic}
  (or a $V$-\emph{statistic})
  for the transformation $T$ and  the kernel $f$. Notice that
  the same class of statistics is determined by
  \emph{symmetric} kernels, so we will assume that $f$ is symmetric
  whenever  it is needed.\\

At first glance the summands in \eqref{sum} can be defined in two
steps. Firstly, the functions $(x_1,\dots,x_d)\mapsto
 f(T^{i_1}x_1,...,T^{i_d}x_d)$ can be obtained using the
dynamics coordinatewise; secondly, they should be restricted to
the main diagonal of $X^d$. The second step, however, requires
some care. Analysis and clarification of the concept of
restriction became another important objective of this work. This
is a crucial point determining substantially the approach in the
present paper. If $f: X^d \to \R$ is a measurable function on the
Cartesian power $(X^d, \mathcal F^{\otimes d},\mu^d)$, it is
viewed, as usual, not as an individual function, but rather as an
equivalence class of individual functions any two of which agree
on some set of measure $1$. Such an equivalence class, in general,
does not have a well-defined restriction to a subset of measure
zero, like the main diagonal is in the case of the atomless space
$(X,\mathcal F,\mu)$. However, some equivalence classes may
contain individual functions with well-defined restrictions (for
example, continuous functions, assuming that $X$ is the unit
interval with the Lebesgue measure $\mu$). A simple but important
observation made in this article is that suitable nice functions
on product probability spaces can be described in purely
measure-theoretical terms. The key concept here is the
\emph{projective tensor product} of Banach spaces. First we show
that, under appropriate assumptions, the elements of a respective
abstract Banach space can be viewed as functions from some
$L_p(\mu^d)$. In particular, every such a function determines an
equivalence class discussed above. Analogously to the situation
with continuous functions, nice representatives (non-unique) can
be found within every such equivalence class; in view of specific
properties of projective tensor products, they can be represented
by absolutely convergent series of the products of functions in
separate variables. Furthermore, such 'special representatives'
can be restricted to the main diagonal in a correct way. Notice,
that the main diagonal is considered here as a probability measure
space whose measure is the image of $\mu$ under the map $x \mapsto
\underbrace{(x, \dots,x\,)}_{d \, \, \text{times}}$\,; correctness
means here that possible uncertainty in the choice of the
restricted function concerns only sets of measure $0$ on the
diagonal. We emphasize that this procedure of 'naive restriction'
applies to 'special representatives' of equivalence classes only.
Different choice of a representative within the same equivalence
class may lead to misunderstandings which can be observed in the
literature. In the present paper, however, another approach to the
restriction problem is developed. Using general properties of
projective tensor products, a \emph{restriction operator} is
defined. We will see that this operator agrees with the 'naive
restriction' in the case of the sums of product functions and
their proper limits. On the other hand, for every equivalence
class of measurable functions discussed above, the restriction
operator can (or can not) be applied to the \emph{entire
equivalence class} and sends it, if applicable, to an equivalence
class of functions on the diagonal; thus, no special choice of a
representative within the class is needed. Moreover, we show in
Proposition \ref{embrestr} that the correct restriction can be
obtained as the result of a natural procedure combining
approximation and regularization (compare with the Steklov
smoothing operators and Theorem 8.4 in \cite{GoKrGo2000}).
Finally, we obtain, along with the \emph{correctness of the
restriction}, its \emph{continuous dependence on the kernel}; this
continuity is critical for our approach. The above discussion
introduces the following result which summarizes Lemma 1 and a
particular case of Proposition 1 in Section 2 where also some
information on projective tensor products can be found. We denote
by $L_p(\mu^d)$ the space $L_p\bigl(X^d, \mathcal F^{\otimes d},
\mu^d \bigr)$ and by
$|\cdot|_p$ the norm in any space $L_p$.\\

\noindent{\bf Statement A:} \emph{Let $p\in [1,\infty)$  and
  $dr=p$. Then the projective tensor product
  $L_{\,p,\,\pi}(\mu^d)$ of $d$ copies of $L_{\,p}(\mu)$ is
  contractively embedded
  into $L_{\,p}\,(\mu^d)$ as a dense subspace. The embedding is given  by a linear map
  sending an elementary tensor
  $f_1\otimes...\otimes f_d$
  to the function $(x_1,\dots,x_d) \mapsto f_1(x_1) \cdots
  f_d(x_d)$.
Moreover, the linear map $D_d$
  defined on elementary tensors $f_1\otimes...\otimes f_d$
  by the relation
  $$D_d(f_1\otimes...\otimes f_d)(x)= f_1(x) \cdots
  f_d(x), \, x \in X,$$
  is a norm $1$ linear map of
  $L_{p,\,\pi}(\mu^d) \subset L_p(\mu^d)$
  to $L_r(\mu)$.
  }\\

We shall see that the map $D_d$ is compatible with the dynamics defined
 by $T$ in the sense that for every $x  \in X$ and $n_1,\dots,n_d \in
   \Z_+$
 $$ D_d \bigl((f_1\circ T^{n_1})\otimes...\otimes
(f_d\circ T^{n_d})\bigr)(x)\!=\! (f_1\circ T^{n_1})(x)\cdots
(f_d\circ T^{n_d})(x).$$ For $\mathbf k=(k_1,...,k_d)$ and
$\mathbf n =(n_1,...,n_d)$ we use the notation $\mathbf k
\boldsymbol{<} \mathbf n $ ($\mathbf k \boldsymbol{\le} \mathbf n
$)
 if $ k_l<n_l$ (respectively, $ k_l \le n_l$) for every
$l=1,...,d$; we set for a function $f: X^d \to \R$ $$\bigl(V^{\bf
k}f\bigr)(x_1, \dots,x_n)=f(T^{k_1}x_1,\dots,T^{k_d}x_d), \,
x_1,\dots,x_d \in X.
$$
The operators $V^{\bf k}$ act on every space $L_{\,p}(\,\mu^d\,)$
and also on every  $L_{p,\,\pi}(\mu^d\,)$ $(1 \le p  \le \infty)$.
So do the (pre-)adjoint operators $V^{*\,{\bf k}}$ (details are
contained in Section 2).\\

 Statement A leads to the following version of the multivariate ergodic
 theorem (Corollary 2 in Section 3).\\

 \noindent{\bf Statement  B:} {\it Let $p=dr$, $1\le
r,p < \infty$. Then for $f \in L_{p,\,\pi}(\mu^d)$ we have
\begin{equation}\label{mises} \frac 1{n_1\,n_2\,...\,n_d}\sum_{\mathbf 0 \boldsymbol{\le}
\mathbf k\boldsymbol{<} \mathbf n}
  D_dV^{\bf k}f \to D_d E^{\otimes d}_{\text{inv},\,\pi}\,f
\end{equation}
almost surely and in the norm of $L_r(\mu)$ as $n_1, \dots, n_d \to \infty.$\\
Here
$E_{\text{inv}}$ is the conditional expectation operator
with respect to the $\s$-algebra of \,$T$-invariant sets, and
$E^{\otimes d}_{\text{inv},\,\pi}$ is
the $d$-th projective tensor power of $E_{\text{inv}}$}. \\

The distributional limit theorems rely on the Hoeffding
decomposition. For every $m \in \{1,\dots,d\}$ let
$L_p^{sym}\,(\mu^m)$  be the subspace of symmetric elements of
$L_p\,(\mu^m)$, $\mathcal S^m_d$ be the collection of all subsets
of $\{1,...,d\}$ of cardinality $m$ and, for every $S \in \mathcal
S^m_d$, let $\pi_S$ be the projection map from $X^d$ onto $X^m$ which
only keeps coordinates with indices in $S$. The \emph{symmetric
Hoeffding decomposition} asserts the existence of operators $R_m:
L_p^{sym}\,(\mu^d)\to L_p^{sym}\,(\mu^m)$ such that every $f \in
L_p^{sym}\,(\mu^d)$ can be represented in a unique way in the form
\begin{equation*}
f=\sum _{m=0}^d \sum_{S\in \,\mathcal S^m_d}(R_m f)\circ \pi_S
\end{equation*}
(see Section \ref{4} for details). The same or analogous notation
will be applied to the spaces $L_{p,\,\pi}(\mu^d)$. \\
In the following Statement C (Theorem 2 in Section 7)) we assume
that $T$ is an \emph{exact} transformation in the sense that
$\bigcap_{n \ge 1} T^{-n} \mathcal F = \mathcal N$, where
$\mathcal N$ is the trivial sub-$\s$-field of $\mathcal F$. Let $
E$ denote the expectation operator. Using the Hoeffding
decomposition and applying to every of its components the
multiparameter martingale-coboundary representation \cite{Go2009}, we prove\\

\noindent{\bf Statement C\,:} {\it Let $T$
be an exact transformation and $f\in
L_2^{sym}(\,\mu^d)$
be a real-valued kernel. Assume that for every $m=1, \dots, d$, $R_m
f \in
 L_{2m,\,\pi\,}^{sym}(\,\mu^m)$ and the series
\begin{equation} \label{thesum2}
\underset{\mathbf 0 \boldsymbol{\le} \, \mathbf k
\boldsymbol{<}\boldsymbol{\infty}}{\sum}\,\,\, V^{*\,\mathbf k}
R_m f \, \biggl(\overset{\mathrm{def}}{=}\, \underset{\mathbf n
\boldsymbol{\to} \boldsymbol{\infty}}{\lim}\,\,\,
\underset{\mathbf 0 \boldsymbol{\le} {\bf k}\boldsymbol{< }\mathbf
n }{\sum}\,\,\, V^{*\,\mathbf k}R_m f \biggr)
\end{equation}
 converges in $
L_{2m,\, \pi}(\,\mu^m)$ (here $\mathbf k= (k_1,\dots,k_m)$,
$\mathbf n=(n_1,\dots,n_m)$).
 Then
$$ V_n^{(d)}f \overset{\mathrm{def}}{=}\,\frac{1}{n^ {d-1/2}}
\underset{0\,\le\,\, k_1,\,\dots,\,k_d \,\le\,\,
n-1}{\,\,\,\sum}\, D_d \, V^{(k_1,\, \dots,\,k_d)}(f-R_0 f)$$
converges in distribution to a centered Gaussian random variable with variance $d^2\sigma^2(f) \ge 0$, where
\begin{equation*}
\sigma^2(f)= \biggl|\sum_{k=0}^{\infty} V^{*\, k} R_1 f
\biggr|_2^2 - \biggl|\sum_{k=1}^{\infty} V^{* k} R_1
f\biggr|_2^2\ge 0.
\end{equation*}
The convergence of the second  moments
\begin{equation*}
E(V_n^{(d)}f)^2\underset{n \to \infty}{\to}{{d\,}^2 \sigma^2(f)}
\end{equation*}
holds as well.}  \\

This Central Limit Theorem (CLT) is complemented by Theorem 3 in
Section 7 which asserts, under weaker assumptions, only the
convergence of the first absolute moments (besides the convergence
to the Gaussian distribution). Last, in Theorem 4 of Section 8, we
prove the following distributional result when $d=2$ and $f$ is a
symmetric \emph{canonical} kernel
(that is $R_0\,f=0$ and $R_1\,f=0$).\\

\noindent{\bf Statement D:} {\it Let $d=2$. For every canonical
$f$ satisfying the assumptions in Statement C there exists an
absolutely summable real  sequence $(\lambda_m)_{m \in \N}$ such
that the random variables
$$ \frac 1 n\sum_{0\le \, k_1\!,\,k_2 \le\,n-1} DV^{(k_1,\,k_2)}f$$
converge in distribution, as $n\to \infty$, to
$$ \xi=\sum_{m=1}^\infty \lambda_m \eta_m^2$$
where $(\eta_m)_{m\ge 1}$ is a sequence of independent standard
Gaussian random variables.
 Moreover,
\begin{equation*}
E \,\Bigl(\, \frac{1}{n}\,  \sum_{0 \le\, i_1\!, \,i_2 \le\, n-1}
D_2 V^{(i_1,\,i_2)} f \Bigr) \underset{n \to \infty}{\to}
\sum_{m=1}^{\infty}\lambda_m .
\end{equation*}}\\

The main limit theorems are presented with proofs in Sections 3, 7
and 8. Section 2 contains necessary preliminary material; in
particular, the restriction operator is introduced there. The
Hoeffding decomposition and filtrations are discussed,
respectively, in Sections 4 and 5. Section 6 contains the main
part of the preparatory work for the rest of the paper. It is here
that the martingale decomposition undergoes the projective tensor
multiplication, leading from the classical Burkholder martingale
inequality  to upper bounds for certain multiparameter sums. These
bounds allow (Section 7) to neglect the influence of higher degree
summands in the Hoeffding decomposition to the asymptotic behavior
when proving the CLT in the non-degenerate case. They are also
applied in Section 8 in the proof of Statement D to show that the
contribution of ``partial coboundaries'' vanishes in the limit; this
reduces the proof to the particular case of a kernel with maximal
possible martingale difference properties. Some examples (in fact,
mostly general results treating entire classes of stationary
processes and kernels) are collected
in Section 9.\\
The above stated results, along with their modification for the
case of invertible transformations (see Remark \ref{adapt}) and
the examples in Section 9, clearly show that a substantial part of
the limit theory for $V$-statistics of stationary processes can be
developed, basing exclusively on projective tensor products and
martingale approximations. The latter is presented only  in its
original primitive form (moreover, only the the adapted case is
considered). Using more recent developments could substantially
relax many assumptions in the paper. Many other limit results can
be established similarly or at the expense of small additional
efforts. However, we believe that this presentation is more
suitable for introducing the subject.

\noindent\begin{remark}
For a given function
$f$ defined on $(X^d, \mathcal F^{\otimes d},\mu^d)$,
a natural question arises to decide
whether $ f \in L_{p,\,\pi}(\mu^d)$
and to bound
its norm.
For $d=2$ and some $p \in (1,
\infty],$ $p^{\,\prime} \in [1,\infty)$,
$p^{-1}+(p^{\,\prime})^{-1}=1$,
an equivalent question is whether the integral  operator from
$L_{p^{\,\prime}}$ to
 $L_{\,p}$ with the kernel $f$ is \emph{nuclear} \cite{Ry2002}.
 There is an extensive literature on the topic, especially
 on nuclear (or \emph{trace class}: see \cite{ReSi1980} and
 also \cite{Ry2002} where Exercise 2.12 shows the difference between the complex
 and the real cases)
operators in Hilbert spaces.
Criteria for
integral operators
to be nuclear
can be traced  back to classical papers of Fredholm and
Carleman (see monographs \cite{GoKr1969,GoKrGo2000} and references
therein; in \cite{GoKrGo2000}  also nuclear operators in Banach
spaces are considered).
A special class consists of positive
semidefinite kernels.
For example, the well-known Mercer's theorem
implies  $f \in L_{2,\,\pi}(\mu^2)$ for such kernels under the
additional assumption that $X$ is a compact space and $f$ is
continuous.\\
To the best of
our knowledge, for $d \ge 3$, much
fewer literature
exists on
this topic. The main tool here is
the expansion of $f$ into a functional series whose summands are
products of sufficiently regular functions in separate variables
$x_1,\dots,x_d$ (see Proposition \ref{suffic} and Section
\ref{applic} for some examples).
\end{remark}
\begin{remark} \label{ustat} The $U$-statistics (that is,
for symmetric kernels $f$, the off-diagonal modification of sums
\eqref{sum}) are mentioned but not treated in the present paper.
Under some strengthening  our assumptions (the series in
\eqref{thesum2}, \eqref{thesum1} and \eqref{conv8} should converge
unconditionally; for example, this will be the case if we are in
the  position to check the assumptions of Proposition
\ref{suffic}) the conclusions of Theorems 2, 3 and 4 can be
reformulated for $U$-statistics. Notice that both advantages of
$U$-statistics compared to $V$-statistics in the i.i.d. case (to
be unbiased estimates of the mean value of the kernel with i.i.d.
arguments; to require weaker assumptions imposed on the kernel) in
general are no longer valid in the dependent case.
\end{remark}
\subsection{Some history and earlier results}
The theory of $U$- and $V$-statistics for i.i.d. variables is well
developed (see \cite{KoBo1994, DeGi1999} and references therein).
Degenerate von Mises statistics for independent variables have
first been treated by von Mises in \cite{vMis1947} and Filippova
in \cite{Fil1962}. Neuhaus \cite{Neu}
  proved a functional form of the weak convergence  for degenerate kernels of degree $2$.
 Although he dealt with the  $U$-statistics only, the method applies as well to von
Mises statistics with properly modified limit distributions.  In
 \cite{DeGrKe85} the functional form of Filippova's result is obtained
 with the distributional limit presented by multiple
 stochastic integrals with respect to the Kiefer--M\"uller process.
Many fine results on  $U$-statistics (maximal inequalities, large
deviations, functional CLT) are included or surveyed in
\cite{DeGi1999} and \cite{Maj2005}.

 For
non-independent random variables some progress has been made for
weakly dependent and associated processes (see \cite{DeTa},
\cite{De} and references therein). More generally, the Strong Law
of Large Numbers (SLLN ) for von Mises statistics of an ergodic
stationary real-valued processes $\xi=(\xi_n)_{n\ge 0}$
with
one-dimensional distribution $\nu$ has been treated in
\cite{Aa96}, where it is shown, among other important results and
interesting examples, that almost surely we have
\begin{equation} \label{misproc}
n^{-d}\!\!\!\!\!\sum_{0 \le i_1<\,n,...,\,0 \le \,i_d\,<n}
F(\xi_{i_1},...,\xi_{i_d})\underset{n \to \infty}{\to}\int_{X^d}
F(x_1, \dots,x_d) \nu(dx_1) \cdots \nu(dx_d),
  \end{equation}
 the assumptions ranging from continuity of
the kernel $F$ to the weak Bernoulli property of $\xi$. One of the
results in \cite{Aa96} on von Mises statistics is a
 SLLN
under the assumption that
 the kernel is  bounded by  a product of functions in separate variables.
  In case of functionals of mixing processes
 a form of the SLLN has been proven in \cite{BoBuDe2001} which is not
 contained in \cite{Aa96}.  In almost all other papers the CLT (sometimes together
 with its functional form) has been considered.
  Yoshihara \cite{Yosh} was the first to give a probabilistic treatment of the
 CLT question when the process is \emph{absolutely
 regular}. Other \emph{mixing conditions} are investigated in \cite{Bab1, Bab1989,
 BoBy2006-1, BoBy2006-2, BoVo2008, DeWe, Kha, Sha, Vol2011, Yosh92}.
Functionals of absolutely
 regular processes have been studied
 in \cite{DeKe84}. In \cite{DeKe86} these results were used to
 construct a new type of asymptotically distribution free  confidence intervals for the
correlation dimension (see
 \cite{GraPro}). Later many limit results have been considerably improved in \cite{BoBuDe2001} and
 \cite{BoBuDe2002} by establishing a functional form of the central limit theorem.
 In the weakly dependent case we mention the works of Babbel
 \cite{Bab1,Bab1989} and Amanov \cite{Ama} where  various types of mixing conditions are considered,
 including strong mixing. The above list is incomplete, more information is contained in the
 surveys \cite{DeTa} and  \cite{De}.

 Notice that in a recent
paper \cite{LeNe2011}, independently of our research, for a
certain class of canonical symmetric kernels of degree 2 (in
\ref{stat}.1 we call them \emph{martingale kernels}) a limit
distribution of $V$-statistics is derived which has  the same form
as in the i.i.d. case. This conclusion agrees with ours in
Statement D above; the result in \cite{LeNe2011} is a rather
particular case of our Statement D (see \ref{stat}.1 for more
details). The paper \cite{LeNe2011} and the subsequent papers
\cite{LeNe2012, Le2012} also develop impressive statistical
applications of this and other limit results; some new, compared
to \cite{LeNe2011}, limit theorems in \cite{LeNe2012, Le2012} are
developed by means of methods different from those used in the
present paper; the corresponding assumptions about the process
include some decay of the Kantorovich distance between the
conditional and the unconditional distributions of the process
given its past; also some form of the Lipschitz condition is
imposed on the kernel. The spectral decomposition of the kernel
or, alternatively, its approximation by Lipschitz continuous
wavelets are used there to derive the results.\footnote{Though all
this creates a favorable environment for employing our Proposition
\ref{suffic}, we do not investigate this possible application in
the present paper.}

\section{Preliminaries} \label{2}

\subsection{Multiparameter actions} \label{2.1}
 Let $T$ be a
measure preserving transformation of a probability space
$(X,\mathcal F, \mu)$ (which is assumed to be standard, that is a
Lebesgue space in the sense of Rokhlin \cite{Ro1961}). For every
$p \in [1,\infty]$ we set $L_p(\mu)= L_p(X,\mathcal F, \mu) $,
choosing $\C$ as the field of scalars and denoting by
$|\cdot|_{\,p}$ the norm of $L_p(\,\mu)$. Define an isometry
$V:L_p(\mu)\to L_p(\mu)$ by the relation $Vf= f\circ T$. For every
$p \in [1,\infty)$ let $V^*:L_{p\,'}(\mu)\to L_{p\,'}(\mu)$ be the
adjoint operator of $V:L_p(\mu)\to L_p(\mu)$ where
$p^{-1}+p\,'^{-1}=1$. The preadjoint operator (acting in
$L_1(X,\mathcal F, \mu) $) of the operator $V: L_{\infty}(\mu)\to
L_{\infty}(\mu)$ will be loosely called  the adjoint of $V$  and
denoted by $V^*$, too, whenever this does not lead to a
misunderstanding. Analogous notations and agreements will be
applied to other measure spaces, their transformations and related
operators.

 For every
$i=1,\ldots,d$ let $(X_i,\mathcal F_i, \mu_i, T_i)$ be a
probability space with a measure preserving transformation $T_i$;
let $V_i, V_i^*$ be the corresponding operators. We assume that
these spaces are copies of $(X,\mathcal F,\mu).$ The direct
product $\prod_{1\le i\le d}(X_i,\mathcal F_i, \mu_i)$ will be
denoted by $(X^d, \mathcal F^{\otimes d}, \mu^d).$ Unlike the
spaces, the transformations $T_1,\dots,T_d$ can be different;
however, from Section \ref{8} on we assume that they are copies of
the same transformation $T$. The notation $L_p(\mu^d)$ should be
understood correspondingly. Let $\Z_+=\{0,1, \ldots\}.$ For every
$\mathbf{n}$ $=(n_1, \ldots,n_d)$ $ \in \Z^d_+$ we set
$T^{\mathbf{n}}(x_1, \ldots,
x_d)=(T_1^{n_1}x_1,\ldots,T_d^{n_d}x_d).$ Define a representation
of the semigroup $\Z^d_+$ by isometries in $L_p(\mu^d)$ via
\begin{equation*}
V^{\mathbf n} f =f \circ T^{\mathbf{n}},\, f \in L_p(\mu^d).
\end{equation*}
We do not assume that the transformation $T$ is invertible. The CLT proved
below will hold for the
 class of essentially noninvertible $T$ (known
as {\em exact transformations}\,).  The family of adjoint
operators $(V^{\mathbf n*})_{\mathbf n \in \Z^d_+}$ is also a
representation of $\Z^d_+$ (by {\em coisometries} in this case).
Note that these two representations do not commute with each other
in the noninvertible case (otherwise they clearly commute).
However, if $\mathbf{e}_1, \dots, \mathbf{e}_d $ denote the
standard basis of $\Z_+^d$, the operators $V^{\mathbf{e}_i}$ and
$V^{*\,\mathbf{e}_j}$ commute for $i \neq j$ because they act on
different coordinates in $X^d$. This will be used in the proof of
Lemma \ref{martcobdec}.
\subsection {Tensor products and products of
functions}\label{2.2} We discuss here conditions on kernels under
which $V$-statistics are well-defined. Recall the concept of the
projective tensor product of Banach spaces \cite{Ry2002,
DeFl1993}. The main field is assumed to be $\C$ or $\R$.

 Let $B_1, \dots,B_d$ be Banach spaces with
norms $|\cdot|_{B_1}, \ldots,|\cdot|_{B_d}$ and let $B_1 \otimes
\dots \otimes B_d$ be their algebraic tensor product. Elements of
$B_1 \otimes \cdots \otimes B_d$, representable in the form
$f_1\otimes  \cdots \otimes f_d$, are called {\em elementary
tensors}. The {\em projective tensor product} of $d \ge 2$ Banach
spaces denoted by $B_1 \hat{\otimes}_{\pi} \cdots
\hat{\otimes}_{\pi} B_d$ is, by definition, the completion of the
algebraic tensor product with respect to the {\em projective norm}
defined as the supremum of all \emph{cross norms} on $B_1 \otimes
\cdots \otimes B_d$.  Recall that a norm on $B_1 \otimes \dots
\otimes B_d$ is said to be a cross norm whenever it equals
$\prod_{i=1}^d |f_i|_{B_i}$ for every elementary tensor
$f_1\otimes  \cdots \otimes f_d.$

 Recall that for every $i=1,\ldots,d$ $(X_i,{\mathcal F}_i,\mu_i)$
 is a copy of $(X,\mathcal F,\mu)$. For $p_1,
\dots, p_d \in [1, \infty]$ we denote by
$|\cdot|_{p_1,\,\ldots,\,p_d,\,\pi}$ the norm of the space
$$ L_{p_1}(X_1,\mathcal F_1, \mu_1)  \hat{\otimes}_{\pi} \cdots
 \hat{\otimes}_{\,\pi}  L_{p_d}(X_d,\mathcal F_d,
\mu_d).$$
 If $p_1= \ldots=p_d =p \in [1,\infty]$,  the above projective
 tensor product and its norm will be denoted by $L_{\,p,\,\pi}(\mu^d)$
 and $|\cdot|_{\,p,\,d,\,\pi}$, respectively.
 We show in the following lemma that  $ L_{\,p,\,\pi}(\mu^d)$
 can be thought of as a subspace of $L_{\,p}(\mu^d)$; hence, its elements
 can be viewed as functions on $X^d$. Some useful properties of these
 functions are established in \ref{2.3}.
\begin{lemma}\label{embedding}
 For
every $p \in [1,\infty]$ there exists a unique linear map
$$J_d: L_{p,\pi}(\mu^d)\to L_p(\mu^d)$$
  of norm $1$ which sends every elementary tensor
$f_1\otimes \cdots \otimes f_d$ to the function $(x_1, \dots,
x_d)\mapsto f_1(x_1) \cdots f_d(x_d)$. Moreover, $J_d$ maps $
L_{p,\,\pi}(\mu^d)$ into $L_p(\mu^d)$ injectively. For $ p \in [1,
\infty)$ $J_d \bigl(L_{p,\,\pi}(\mu^d)\bigr)$ is dense in
$L_p\,(\mu^d)$.
\end{lemma}
\begin{proof} The case $d=1$ is trivial, so we assume $d \ge 2$. For every
$p\in [1,\infty]$, let us define a linear map $J_d$ of norm $1$,
$$J_d: L_{p\,,\,\pi}(\mu^d)\to L_p(\,\mu^d).$$
When we need to specify $p$ we shall use the notation $J_{\,d,\,p}$.
First, sending every elementary tensor $f_1\otimes \cdots \otimes
f_d$ to the function $(x_1, \dots, x_d)\mapsto f_1(x_1) \cdots
f_d(x_d) ,$ we define a $d$-linear map of norm $1$ from $
L_{\,p}(X_1,\mathcal F_1, \mu_1) \times \cdots \times
L_{\,p}(X_d,\mathcal F_d,\mu_d)$ to $L_{\,p}(\mu^d).$ Then, by a
general property of the projective tensor product
 (see \cite{Ry2002}, Theorem 2.9, for $d=2$; use induction and associativity
 for $d>2$) this map extends to
 $ L_{\,p,\,\pi}(\mu^d) $ uniquely with norm $1$. Denote this resulting map by
 $J_d$.  Its
 image is dense in $L_{\,p,\,\pi}(\mu^d) $ for $p < \infty$ since so is the
 image under $J_d$ of the algebraic tensor product.\\ We
 prove now that $J_d \,(\,=J_{d,\,p})$ is injective. For $p=1$ it
 is so because
 $J_{d,1}$ is an isometric isomorphism between its domain and its range
 (\cite{Ry2002}, Exercise 2.8).
 Let now for some $p > 1$  $I_{1,\,p}: L_p(\mu) \to L_1(\mu)$ and
 $I_{d,\,p}: L_p(\mu^d) \to L_1(\mu^d)$ be the inclusion operators (of norm $1$ each).
 By the \emph{metric mapping property} (\cite{DeFl1993},\,12.1) of the
 projective tensor norm, the inclusion $I_{1,\,p}$  gives rise
 to the norm $1$ mapping $A_d: L_{p,\,\pi}(\mu^d)\to
 L_{1,\,\pi}(\mu^d)$ (notice that $L_{1,\,\pi}(\mu^d)$ and $L_1(\mu^d)$ are
 identified by $J_{d,\,1}$).
  Since the spaces $L_p$ have the \emph{approximation
 property}, the operator $A_d$ is injective as a projective tensor product
 of injective operators $I_{1,\,p}$ (see Corollary 4\,(1), subsection 5.8,
 in \cite{DeFl1993}; then use induction).
 Starting with algebraic tensor products and
passing,
 in view of boundedness of all operators involved, to the
 completions with respect to corresponding norms, we obtain that the mappings
 $J_{d,\,1}\, A_d: L_{p,\pi}(\mu^d) \to L_1(\mu^d)$ and
 $I_{d,\,p}\, J_{d,\,p}: L_{\,p,\,\pi}(\mu^d) \to L_1(\,\mu^d)$ agree.
  Since
 $A_d$ and $J_{\,d,\,1}$ are injective, so is $J_{\,d,\,p}.$
\end{proof}
In view of the properties of $J_d$ we shall, when possible, omit
the symbol $J_d$ and consider  $L_{\,p,\,\pi}(\mu^d)$ as a
subspace of $L_{\,p\,}(\mu^d)$.
 Set for $\mathbf n=(n_1,
\dots,n_d)$
\begin{equation} \label{equivalence}
 V^{\mathbf
n}_{\pi}=V_1^{n_1}{\otimes}_{\pi}\cdots{\otimes}_{\pi}
 V_d^{n_d},  \,\,\,\,\,\, V^{*\mathbf
n}_{\pi}=V_1^{*n_1}{\otimes}_{\pi}\cdots{\otimes}_{\pi}
 V_d^{*n_d}.
\end{equation}
The operators $ ( V^{\mathbf n}_{\pi},  V^{*\mathbf
n}_{\pi})_{\mathbf{n}\in
   \Z^d_+}$ have  properties very similar to those of $ ( V^{\mathbf n},
V^{*\mathbf n})_{\mathbf{n}\in \Z^d_+}$; in particular, they have
norm $1$ with respect to the projective tensor norm. The relations
$J_{d,\,p} V^{\mathbf n}_{\pi}= V^{\mathbf n} J_{d,\,p}$,
$J_{d,\,p} V^{*\mathbf n}_{\pi}= V^{*\mathbf n} J_{d,\,p}$\,,
$\mathbf{n}\in \Z^d_+, $ are obvious for elementary tensors and
immediately extend to the general case. It follows from these
relations that the space $L_{\,p,\,\pi}(\mu^d)$ is preserved by
the operators $ ( V^{\mathbf n}, V^{*\mathbf n})_{\mathbf{n}\in
\Z^d_+}$.  From now on we shall use the notation $ ( V^{\mathbf
n}, V^{*\mathbf n})_{\mathbf{n}\in \Z^d_+}$ also to denote the
restrictions of these families to the space $ L_{\,p,\,\pi}(\mu^d)
\subset L_{\,p\,}(\mu^d)$.
\begin{remark} \label{traceclass} The space $ L_{2,\,\pi}(\mu^2)$
can be identified with the space of  nuclear (or \emph{trace
class})  operators from ${L_2(\mu)}^*$ to $L_2(\mu)$
(\cite{Ry2002}). The operator $J_2$ in Lemma \ref{embedding}
transforms such (integral) operators to their kernels which form a
subspace of $ L_{2}(\mu^2)$.
\end{remark}
\subsection{Restriction to the diagonal} \label{2.3}
  In the following Proposition \ref{restriction}, for  every $p_1,\dots,p_d \in [1, \infty]$ with
 $p_1^{-1}+\dots + p_d^{-1}=1$ and for every
 $f \in L_{p_1}(\mu)  \hat{\otimes}_{\pi} \cdots
 \hat{\otimes}_{\pi}  L_{p_d}(\mu),$  we  define
a function
 $D_df\in L_1(\mu).$ In the case of $1 \le p_1= \dots =p_d=p \le \infty$
 the embedding $J_d$  (Lemma \ref{embedding}) allows us to consider the space
 $L_{\,p}(\mu)  \hat{\otimes}_{\,\pi} \cdots  \hat{\otimes}_{\,\pi}  L_{\,p}(\mu)$
 as a subspace of the $L_{\,p}(\mu^d)$ and interpret its elements as
 functions defined on $X^d$. Then $D_d f$ plays the role of the restriction
of $f$ to the principal
 diagonal $\{(x_1,\dots,x_d): x_1=\dots=x_d)\}\subset X^d$. In this
particular case the term 'restriction'  can be
 justified by an approximation procedure described in Proposition \ref{embrestr}
 below.
\begin{proposition} \label{restriction} Let $p_1, \dots, p_d \in [1, \infty],$ $r \in [1,\infty]$
sa\-tis\-fy
\begin{equation*}
\sum_{i=1}^d \frac{1}{p_i}= \frac{1}{r}.
\end{equation*}
Then
\begin{enumerate}
\item the map $\mathcal D$, sending
 every $d$-tuple $(f_1,\dots,f_d) \in L_{p_1}(\mu)\times \dots$ $\times L_{p_d}(\mu)$
to the function
\begin{equation*}
x \mapsto f_1(x) \cdots f_d(x),
\end{equation*}
is a norm $1$ $d$-linear map from $L_{p_1}(\mu)\times \dots \times
L_{p_d}(\mu)$ to $L_r(\mu);$ \item there exists a unique linear
map (of norm $1$)
\begin{equation*}
D_d: L_{p_1}(\mu)  \hat{\otimes}_{\pi} \cdots  \hat{\otimes}_{\pi}
L_{p_d}(\mu) \to L_r(\mu)
\end{equation*}
such that for every $d-$tuple $(f_1,\dots,f_d) \in L_{p_1}(\mu)\times
\dots \times L_{p_d}(\mu)$
$$D_d(f_1 \otimes \dots \otimes f_d) = \mathcal D (f_1, \dots, f_d).$$
\end{enumerate}
\end{proposition}

\begin{proof}
The first assertion is a consequence of the multiple H\"older
inequality (Exercise 6.11.2 in \cite{DuSchw1958}).  The second one
follows from the linearization property of the projective tensor
products with respect to polylinear maps.  For the case of
bilinear maps see Theorem 2.9 in \cite{Ry2002}; for $d>2$ use
induction and associativity.
\end{proof}

 If $p_1=\dots=p_{\,d}=p,$ the space
 $ L_{p,\, \pi}(\mu^d)= L_{p}(\mu)  \hat{\otimes}_{\pi} \cdots
 \hat{\otimes}_{\pi} L_{p}(\mu)$ is embedded into $L_p(\mu^d)$ by
 the operator $J_d$ (Lemma \ref{embedding}); we omit $J_d$ and treat
 an $f \in L_{p,\, \pi}(\mu^d)$ as a function. For every finite
measurable partition $\mathcal A=\{A_1 \dots, A_m \} $ let us
denote by $\mathcal F_{\mathcal A}$
 the $\sigma-$field of all
possible unions of atoms of $\mathcal A$ and by $E(\cdot|\,
\mathcal A)$ the corresponding conditional expectation.
 Let  $(\mathcal A_n)_{n\ge 1}$ be a refining sequence of
finite measurable partitions $\mathcal A_n=\{A_{1,\,n},$ $\dots,
A_{m_n\!,\,n} \} $ such that $\mathcal F$ is the smallest
$\sigma$-field containing all $\mathcal F_{\mathcal A_n}, n \ge
1.$ Let $I_A$ denote the indicator of the set $A.$

\begin{proposition} \label{embrestr}Let $d\ge 1,$
$p_1=\dots=p_{\,d}=p \in\! [\,d,\infty)$ and $r=p/d.$ Define the
sequence $ (D_{d,\, n})_{n \ge 1}$ of operators $D_{d,\,n}:
{L}_{p,\,\pi}(\mu^d)\to L_r(\m)$ by
$$D_{d,\,n}f =\sum_{i=1}^{m_n} \frac{I_{A_{i,\,n}}}{\mu(A_{i,\,n})^{d}}
\int_{A_{i,\,n}^d}f(x_1,\dots, x_d) \mu(dx_1) \cdots \mu(dx_d).$$
Then $ D_{d,\,n}\underset{n\to \infty}{\to} D_d $ in the strong
operator topology.
\end{proposition}

\begin{proof}
First let us verify (again using the H\"older inequality) that
$D_{d,\,n}$ as a map from $ L_{\,p,\,\pi}(\mu^d)$ to $L_r(\mu)$
does not increase the norms of elementary tensors. From the
relation
\begin{equation} \label{prod}
\begin{split}
D_{d,\,n}(f_1 \otimes \dots & \otimes f_d) =D_d(E(f_1|\mathcal
A_n)\otimes \cdots\otimes E(f_d|\mathcal A_n))\\
&= E(f_1|\mathcal A_n) \cdots E(f_d|\mathcal A_n)
\end{split}
\end{equation}
 it follows that
\begin{equation*}
\begin{split}
&|D_{d,n}(f_1 \otimes \dots \otimes f_d)|_r =|E(f_1|\mathcal A_n) \cdots E(f_d|\mathcal A_n)|_r \\
&\le |E(f_1|\mathcal A_n)|_p \cdots |E(f_d|\mathcal A_n)|_p\le
|f_1|_p \cdots |f_d|_p.
\end{split}
\end{equation*}

By the properties of the projective norm, this implies that the
norm of every $D_{d,\,n}: L_{p,\, \pi}(\mu^d) \to L_r(\m)$ is also
bounded by $1$.

Now, using \eqref{prod}, standard properties of conditional
expectations and the H\"older inequality, we obtain
\begin{equation*}
|D_{d,\,n}(f_1\otimes \cdots \otimes f_d)-f_1 \cdots f_d|_r \le
\sum_{k=1}^d |E(f_k|\mathcal A_n)-f_k|_p \prod^d_{m=1,\,m\neq k}
|f_m|_p.
\end{equation*}
 From the martingale
convergence theorem for the space $L_p$ we conclude that every
sequence $(D_{d,\,n}(f_1\otimes \cdots \otimes f_d))_{n \ge1}$
converges in the norm of $L_r(\m)$ to the function $  f_1(\cdot)
\cdots f_d(\cdot). $ The analogous conclusion holds for finite
linear combinations of elementary tensors. Since the norms of the
operators $ D_{d,\,n}$ are uniformly bounded, the proposition
follows.
\end{proof}
The following corollary will be used in the proof of
 Proposition \ref{ergthm}.
 \begin{corollary} \label{posit} The restriction operator $D_d$ preserves positivity of real
 valued functions.
 \end{corollary}
 Thus, the function $D_d f \in L_r(\mu)$ is a well-defined
substitute for the naive restriction of $f$ to the principal
diagonal. For example, for $ \mathbf n \! = \! (n_1,\dots,n_d)$
the function $D_d V^{\mathbf n}f$ can be viewed as a substitute
for the function $ x \mapsto f(T^{n_1}x, \dots,T^{n_d}x).$

\section {Strong law of large numbers} \label{3}

 \subsection{A multivariate ergodic theorem} \label{3.1}

If $T$ is an ergodic transformation of a probability space, a von
Mises statistic may be considered as an estimate for the multiple
integral of the kernel with respect to the invariant measure.
Consistency is one of the desirable statistical properties of a
sequence of estimates; this raises the question of an appropriate
ergodic theorem.
 Proposition \ref{ergthm}, the main result of this subsection, states such a
theorem in a general setting. It asserts, in the ergodic case, the
convergence of multiparameter sums (\ref{eq:prop3}) to the average
of the kernel with respect to the product measure. This reminds  of a
Wiener-type ergodic theorem (\cite{DuSchw1958}, Theorem 8.6.9)
specialized to the case of $d$ one-parameter coordinatewise
actions on the product of $d$ probability spaces. However, not
only the assumptions, but also the conclusions in these results
are different: unlike the Wiener theorem, our result asserts the
convergence for almost all initial points with respect to a
probability measure which is in general \emph{neither absolutely
continuous} with respect to the product measure (being supported
on the main diagonal) \emph{nor invariant} under the
multiparameter action.

We do not assume here symmetry of the kernel and perform summation
over rectangular coordinate domains (which is common in the
multiparameter ergodic theorems, see \cite{DuSchw1958}, Chapter 8)
rather than over coordinate cubes involved in the definition of
$V$-statistics. In this subsection we consider several
 possibly different $\mu$--preserving transformations $T_1, \dots, T_d$ of the space
 $(X,\mathcal F, \mu)$, using the notation \linebreak
 $T^{(n_1,\dots,n_d)}(x_1,\dots,x_d)\!
 =\!(T_1^{n_1}x_1, \dots,T_d^{n_d}x_d)$
 and $V^{(n_1,\dots,n_d)}f\! =$$ f \circ T^{(n_1,\dots,n_d)}.$  \\
Transformations considered in this subsection in general are not
ergodic, so we need some notations to
include the non-ergodic case.
Recall that $A \in \mathcal F$ is said to be $T-$\emph{invariant}
if $T^{-1}A=A$. For every $l \in \{ 1,\dots,d \}$ let $\mathcal
F_{inv,\,l}$ denote the $\s-$field of all $T_l-$invariant
measurable sets in $(X, \mathcal F, \mu)$, and let $E_{inv,\,l}$ be the
corresponding conditional expectation considered as an operator in
$L_{p_{\,l}}(X, \mathcal F, \mu).$
% Set
%$$E_{inv,\,\pi}= E_{inv,\,1} {\otimes}_{\,\pi} \cdots {\otimes}_{\,\pi} %E_{inv,\,d}.$$
\begin{proposition} \label{ergthm}
Let $p=rd$ for some integer $d\ge 1$  and a real number $r \in
[1,\infty)$.
 Let $T_1,\dots,T_d$ be  measure preserving
transformations of a probability space $(X, \mathcal F, \mu)$ and
$f \in \! L_{p,\pi}(\mu^d)$. Then, with  $\mathbf
n=(n_1,\dots,n_d)$, we have
\begin{equation}\label{eq:prop3} \frac{1}{n_1 \cdots n_d}
\sum_{\mathbf 0 \boldsymbol{\le}\,\mathbf k \boldsymbol{<} \mathbf
n} D_d V^{{\bf k}}f  \underset{n_1,\,\dots,\, n_d\,\to\infty}{\to}
D_d (E_{inv,1} \otimes_{\pi} \cdots
 \otimes_{\pi} E_{inv,\,d})f
\end{equation}
 with probability 1 and in $L_r(\mu)$.
\end{proposition}

\begin{remark} \label{clarlimit} The main point of Proposition \ref{ergthm}
is the convergence with probability 1 in the case $d \ge 2$. As to
the convergence in $L_r$, it is not hard to prove, for every $d
\ge 2$ and $p_1, \dots, p_{\,d}, r \in (1,\infty)$, satisfying
$\sum_{i=1}^d p_i^{-1}= r^{-1}$, the following multiple
statistical ergodic theorem:
\[\frac{1}{n_1 \cdots n_d} \sum_{\mathbf 0 \boldsymbol{\le}\,\mathbf k \boldsymbol{<} \mathbf n}
V^{\mathbf k} \underset{n_1,\,\dots,\, n_d\,\to\infty}{\to}
(E_{inv,1} \otimes_{\pi} \cdots
 \otimes_{\pi} E_{inv,\,d}),\] asserting the
strong convergence in the space
 $ L_{p_1}(X_1,\mathcal F_1, \mu_1)  \hat{\otimes}_{\pi} \cdots$
 $\hat{\otimes}_{\pi}  L_{p_d}(X_d,\mathcal F_d, \mu_d).$ Applying
the operator $D_d$ to the both sides of this relation, we obtain
the convergence in the $L_r$-norm. Choosing $p_1=\cdots=p_d=rd=p$,
the convergence in $L_r(\mu)$ in
Proposition \ref{ergthm}
follows for $d \ge 2$. The proof of Proposition \ref{ergthm} contains a second argument of this fact.
\end{remark}
 The next
 lemma will be used in the proof of Proposition \ref{ergthm}.
\begin{lemma}\label{upperbound} Let $d,p,r$ and
the transformations $T_1,\dots,T_d$ satisfy the conditions of
Proposition \ref{ergthm}. Let, moreover, $p>1$. Then there exists
a constant $C=C(r,d)$ such that for every $f \in L_{p,\pi}(\mu^d)$
we have the inequality
\begin{equation*}
\Biggl| \underset{\substack{1\le n_1 <\infty\\ \dots\\ 1 \le n_d <
\infty}}{\sup} \biggl|\frac{\sum_{k_1=0}^{n_1-1} \cdots
\sum_{k_d=0}^{n_d-1} D_d V^{(k_1,\,\dots,\,k_d)}f}{n_1\cdots n_d}
\biggr| \Biggr|_{\,r} \le C |f|_{p,\,d,\,\pi}.
\end{equation*}
\end{lemma}

\begin{proof} For the proof we will use the bound in \cite{DuSchw1958}, Theorem
8.6.8. Note that this result is the lemma for $d=1$. \\
Let now $d \ge 2.$ According to one of the properties of the
projective tensor norm (\cite{Ry2002}, Proposition 2.8), for every
$f \in L_{p,\,\pi}(\,\mu^d)$ and $\e > 0$ there exists a bounded
family of functions $f_{i,\,l} \in L_{\,p}(\,\mu)$ $(1\le i<
\infty, 1 \le l \le d)$ such that $ f =\sum_{i} f_{i,\,1}
{\otimes}_{\pi} \cdots {\otimes}_{\pi}f_{i,\,d} $ and
\begin{equation*}
\sum_{i} |\,f_{i,\,1}|_p\cdots |\,f_{i,\,d}|_p \le
|\,f|_{p,\,d,\,\pi} + \e.
\end{equation*}
Then we have, using Corollary \ref{posit}, that
\begin{equation*} \label{bound1}
\begin{split}
&\Biggl|\underset{\substack{1\le n_1 <\infty\\ \dots\\ 1 \le n_d <
\infty}}{\sup}\biggl|(n_1\cdots n_d)^{-1} \sum_{k_1=0}^{n_1-1}
\cdots
\sum_{k_d=0}^{n_d-1}D_d V^{(k_1,\,\dots,\,k_d)} f \biggr| \Biggr|_r \\
\le & \sum_{i}\Biggl|D_d\Biggl( \underset{1\le n_1<\infty}{\sup}
\frac{| \sum_{k_1=0}^{n_1-1} V_1^{k_1} f_{i,1}|}{n_1}\cdots
 \underset{1\le n_d <\infty}{\sup} \frac{|
\sum_{k_1=0}^{n_d-1} V_d^{k_d} f_{i,d}|}{n_d}  \Biggr) \Biggr|_r \\
\le & \sum_{i}C|\,f_{i,1}|_{p}\cdots  |\,f_{i,\,d}|_{p} \le
C(|\,f|_{p,\,d,\,\pi} + \e).
\end{split}
\end{equation*}
In the above formulas $V_1, \dots, V_d$ are the dynamical
operators associated with the transformations $T_1,\dots,T_d.$
\end{proof}

\begin{proof}[Proof of Proposition \ref{ergthm}] For $d=1$ the assertions
of the proposition are the classical individual and statistical
ergodic theorems. Let now $d \ge 2$, hence $p \ge 2$. In view of
Lemma \ref{upperbound}, the proof is straightforward. First we
prove the assertions of the proposition for elementary tensors $f
=f_1\otimes \dots \otimes f_d$ with $f_l \in L_{\,p\,}(\mu)$, $1
\le l \le d.$ Then the corresponding normalized $V$-statistic can
be written in the product form
\[ \frac{
\sum_{k_1=0}^{n_1-1} V_1^{k_1} f_1}{n_1}\cdots
  \frac{
\sum_{k_1=0}^{n_d-1} V_d^{k_d} f_d}{n_d}, \] where by the
individual ergodic theorem the $l$-th term in the product
converges to $E_{inv,\,l}\,f_l$ with probability 1. Hence, the
product tends with probability 1 to
$$(E_{inv,\,1}f_1) \cdots (E_{inv,\,d}f_d) =D_d E_{inv,\,1} {\otimes}_{\,\pi} \cdots {\otimes}_{\,\pi} E_{inv,\,d}f.$$
The same conclusion holds for finite sums of elementary tensors
which are dense in the space $L_{p\,,\,\pi}(\mu^d)$. Let now $f
\in L_{p\,,\,\pi}(\mu^d)$. Fix an $\e>0.$ There exists an element
$f_{\e} \in L_{p\,,\,\pi}(\mu^d)$ with $|f-f_{\e} |_{p,\,d,\,\pi}
< \e$ such that the a.s. assertion of the proposition holds for
$f_{\e}$ and with probability 1
\begin{equation*} \label{proof}
\begin{split}
&0 \le \xi \eqdef\underset{\substack{n_1 \to \infty
\\ \dots \\ n_d \to \infty}}{\overline{\lim}}\biggl|\frac{1}{n_1 \cdots n_d}
  \sum_{\mathbf 0
\boldsymbol{\le} \mathbf k \boldsymbol{<}\mathbf n} D_d
V^{{\mathbf k}}f - D_d (E_{inv,\,1} \otimes_{\pi} \cdots
 \otimes_{\pi} E_{inv,\,d})f\biggr|\\
\le & \underset{\substack{n_1 \to \infty \\ \dots \\ n_d \to
\infty}}{\overline{\lim}} \biggl|\frac{1}{n_1 \cdots n_d}
\sum_{\mathbf 0 \boldsymbol{\le} \mathbf k \boldsymbol{<} \mathbf
n} D_d V^{{\bf k}}(f - f_{\e})\biggr|+ \biggl| D_d (E_{inv,\,1}
\otimes_{\pi} \cdots
\otimes_{\pi} E_{inv,\,d})(f_{\e}-f)\biggr|\\
+&\underset{\substack{n_1 \to \infty
\\ \dots \\ n_d \to \infty}}{\overline{\lim}}\biggl|\frac{1}{n_1 \cdots n_d}
\sum_{\mathbf 0 \boldsymbol{\le} \mathbf k \boldsymbol{<}\mathbf
n} D_d V^{{\mathbf k}}f_{\e} - D_d (E_{inv,\,1} \otimes_{\pi}
\cdots
 \otimes_{\pi} E_{inv,\,d})f_{\e}\biggr|\\
&\eqdef \xi_{1,\,\e}+\xi_{2,\,\e}+\xi_{3,\,\e}.
\end{split}
\end{equation*}
Since the operators  $D_d$ and 
 $D_d (E_{inv,1}
\otimes_{\pi} \cdots  \otimes_{\pi} E_{inv,d})$ are of norm $1$, we
have  $|\,\xi_{2,\, \e}|_{\,r} \le \e,$ and, in
view of the individual ergodic theorem and Lemma \ref{upperbound}, $\xi_{3,\,\e}=0$ and  $|\xi_{1,\,\e}|_{\,r} \le C \e.$
This implies $\xi = 0$ which proves the convergence with
probability 1. To establish the $L_r$-convergence, we observe that
we have the convergence with probability 1 along with the
domination by an $L_r$-function given by Lemma \ref{upperbound}.
Hence, we can apply Theorem 3.3.7 in \cite{DuSchw1958}.
\end{proof}

 \subsection{Applications to the SLLN for von Mises statistics}
 We return here to the  assumption that the transformations $T_1, \dots, T_d$ are copies of
 the same transformation $T.$ For simplicity we assume that $T$
is ergodic. Symmetry of the kernel is not  assumed.

\begin{theorem} \label{slln} Let $r=p/d$ for some integer $d \ge 2$ and a real number $p
\ge d.$ Let $T$ be an ergodic measure preserving transformation of
a probability  space $(X, \mathcal F, \mu).$ Assume also that $f
\in  L_{p, \pi}(\mu^d).$ Then, as $n \to \infty,$ the sequence
\begin{equation} \label{expression}
\frac{1}{n^d} \sum_{0\, \le\, k_1,\dots, \,k_d\,\! \le\, n-1} D_d
V^{(k_1,\,\dots,\,k_d)}f
\end{equation}
converges  with probability 1 and in $L_r(\mu)$ to the limit
$$\int_{X^d}(J_df)(x_1, \dots, x_d) \mu (dx_1)\cdots \mu(x_d).
$$
Here $J_d: L_{p,\, \pi}(\mu^d) \to L_p(\mu^d)$ is the operator
introduced in Lemma \ref{embedding}.
\end{theorem}

\begin{proof} The theorem follows from Proposition
 \ref{ergthm}. We only need to identify the limits. Since
 the limit expressions given in Proposition \ref{ergthm} and in the
 theorem are both continuous in the projective norm, it suffices to
 check that these expressions agree for elementary tensors $f_1 \otimes \dots \otimes f_d$. It is straightforward to
 check that in the ergodic case both expressions reduce to
 $Ef_1 \cdots Ef_d$,
 where $E$ denotes the integral with respect to $\mu.$
\end{proof}

\begin{corollary} \label{particular}
In the case  $p=d$ Theorem \ref{slln} applies and gives the
convergence with probability 1 and in $L_1(\mu).$
\end{corollary}
\begin{remark} \label{order} Examples show that it is possible to extend
the class of kernels to which the conclusion in Corollary
\ref{particular} applies to such kernels $f\in L_{p}(\,\mu^d)$
which can be "sandwiched" between decreasing and increasing
sequences of some $ L_{p, \pi}(\,\mu^d)$-kernels whose common $
L_{p}(\,\mu^d)-$limit is $f$ (notice that bounding by products
plays some role in \cite{Aa96}). This indicates that probably more
appropriate functional spaces can be found in order  to treat the
SLLN.
\end{remark}

\begin{corollary} \label{sllnser}
Let $T$ be an ergodic measure preserving transformation of a
probability  space $(X, \mathcal F, \mu)$ and let
$(e_k)_{k=0}^{\infty}$ be a sequence of functions in $ L_d\,(\mu)$
such that $e_0\equiv 1$ and for every $k \ge 1$ $|\,e_k|_{\,d}=1$,
$\int_X e_k(x) \mu(dx)=0$. Let
 $f \in L_d(\,\mu^d)$ admit the representation
\begin{equation*} \label{}
f(x_1,\dots,x_d)=
\sum_{\mathbf 0 \boldsymbol{<} \mathbf k
\boldsymbol{<}\boldsymbol{\infty}}
\l_{\bold k}(f)\, e_{k_1}(x_1)\cdots e_{k_d}(x_d)
\end{equation*}
for some family $(\l_{\mathbf k}(f))_{\mathbf 0
\boldsymbol{<}\mathbf k \boldsymbol{<}\boldsymbol{\infty}}$
satisfying the condition \begin{equation*} \label{sumcoeff}
\sum_{\mathbf 0 \boldsymbol{<}\mathbf k
\boldsymbol{<}\boldsymbol{\infty}} |\,\l_{\,\mathbf k}\,(f)\,|\, <
\infty.
\end{equation*}
Then Corollary \ref{particular} applies to $f$.
\end{corollary}
\begin{proof}
The series representing $f$ obviously converges in $ L_{p\,,
\,\pi}(\mu^d) $, and the corollary follows.
\end{proof}
\section{The Hoeffding decomposition} \label{4}
In this section we recall well-known properties of the Hoeffding
decomposition for kernels  in the spaces $L_{p\,}$, omitting
proofs (see \cite{DyMa1983}  for the proofs in the symmetric
case). It is not hard to see that the results and formulas related
to this decomposition (both general and symmetric) apply also to
the spaces $ {L}_{p\,,\,\pi}$
 and, in case $\mu_1=\cdots= \mu_d=\mu$, to their symmetric
 subspaces.
\subsection{The Hoeffding decomposition for general kernels} \label{4.1}
 Let \linebreak $(X_1,\mathcal F_1, \mu_1), \dots,$ $(X_d,\mathcal F_d,
\mu_d)$ be  probability
 spaces. We do not assume in this subsection that all
$(X_l,\mathcal F_l, \mu_l),$ $ l=1,\dots,d,$ are copies of the
same probability space. Let $\mathcal S_d$ ($\mathcal S^m_d$) be
the set of all subsets (respectively, of all $m$-subsets) of $\{1,
\dots, d\}.$ For every $S \subset \{1, \dots, d\}$ we define
$$
(X^S, \mathcal F^{\otimes S}, \mu^S) = \biggl(\,\prod_{l \in
S}X_l, \bigotimes_{l \in S}\!\mathcal F_l, \prod_{l \in
S}\mu_l\biggr),\, L_p(\mu^S) = L_{p\,}(X^S, \mathcal F^{\otimes
S}, \mu^S).
$$
Denoting  the conditional expectation with
respect to a $\s-$field $\mathcal G \subset \mathcal F$ by $E^{\mathcal G}$
and   the projection map from $X^{\{1,\dots, d\}}$ onto
$X^{\{l\}}=X_l$ $(l=1, \dots, d)$ by
$\pi_l$, we set for every $S \in
\mathcal S_d$
\begin{equation*}
\mathcal F^S=\bigvee_{l \in S} \pi_l ^{-1}(\mathcal F_l), \, \,
E^S=E^{\mathcal F_S}, \,\,  \check{E}^l=E^{\{1, \dots,d\}
\setminus \{l\}}\, .
\end{equation*}
In other terms, applying $\check{E}^l$, one integrates out the
$l-$th
variable. \\
The identity operator $I$ in $L_{p}(\mu^{\{1...d\}})$ $(p \in [1,
\infty])$ decomposes as
\begin{equation*}
 I=\prod_{l=1}^d \bigl(\check{E}^{\,l}+(I-\check{E}^{\,l})\bigr)
=\sum_{m=0}^d \sum_{S \in \mathcal S^m_d} Q_S,
 \end{equation*}
where $Q_S=  \prod_{l \notin S} \check{E}^{\,l} \prod_{l' \in
S}(I-\check{E}^{\,l'})$. In general, the Hoeffding decomposition
assigns to every $f \in
 L_{p}(\mu^{\{1,...,\,d\}})$ the family $(R_S
f)_{S \in \mathcal S_d} $  such that
\begin{itemize}
\item[i)] for every $S \in \mathcal S_d$ \
 $R_S f \in L_p(\mu^S);$
\item[ii)]for every $S=\{l_1, \dots, l_m\} \in \mathcal S^m_d$
\begin{equation*}
(R_S f) \circ \pi_S = Q_S f,
\end{equation*}
where $\pi_S:X^d \mapsto  X^S$ is defined by
$\pi_S(x_1,\dots,x_d)=(x_{l_1}, \dots, x_{l_m});$ \item[iii)]
every $R_S f$ is \emph{canonical} (or, using an alternate
terminology, \emph{totally degenerate}) that is for every $l \in
S$, $f \in L_p^{\{1,\dots,\,d\}}$
\begin{equation*}
\check{E\,}^l\bigr((R_Sf)\circ \pi_S\bigr)=0.
\end{equation*}
\end{itemize}
Kernels of the form $(R_Sf)\circ \pi_S$ will be also called
canonical. $R_Sf$ (or $(R_Sf)\circ \pi_S$) is said to have the
\emph{degree} $m$ whenever it  does not vanish identically and $S
\in \mathcal S^m_d$.  Every kernel $f \in
 L_{p\,}(\mu^{\{1, \dots,d\}})$ can be represented in a unique way
 as a sum of canonical kernels (the Hoeffding decomposition) as
 follows
\begin{equation}\label{split}
f=\sum _{m=0}^d \sum_{S\in\,\mathcal S^m_d}(R_S f)\circ \pi_S.
\end{equation}
As  said before, the Hoeffding decomposition also holds for $
L_{p, \pi}(\,\mu^{\{1, \dots,\,d\}})$
$\overset{\mathrm{def}}{=}L_{p}(\,\mu_1)$ $\hat{\otimes}_{\pi}
\cdots$ $\hat{\otimes}_{\pi}L_{p}(\,\mu_d),$ and we shall use the
above notation
for the operators on these spaces as well.\\

The \emph{degree} of a kernel $f$ with decomposition \eqref{split}
(or the decomposition \eqref{symsplit} below) is, by definition,
the smallest degree of non-vanishing summands in
\eqref{split}. A kernel $f$ in \eqref{split} is called
\emph{degenerate} if the degree of $f-R_{\emptyset}f$ is greater
than $1$ and \emph{non-degenerate} if it equals $1$.

\subsection{The Hoeffding decomposition of symmetric kernels} \label{4.2}
 We assume
in this subsection that all spaces  $(X_l,\mathcal F_l,
\mu_l),$ $l=1, \dots,d,$ are copies of the same probability
space  $(X,\mathcal F, \mu).$ $L_{p\,}(\,\mu^d\,)$ and $
L_{p\,,\,\pi\,}(\,\mu^d)$ denote, respectively, the usual
$L_p$--spaces of the product of $d$ identical probability spaces
and the projective tensor product $
\underbrace{L_p(\mu)\hat{\otimes}_{\pi} \cdots
\hat{\otimes}_{\pi}L_p(\mu)}_{d \, \, \text{times}}$ with the
norms $|\cdot|_p$ and $|\cdot|_{p,\,d,\,\pi},$ respectively. There
is an isometric action of the symmetric group $\mathbf S_d$ by
permutations of the multipliers on every of these spaces. The
fixed points of these actions form closed subspaces called
\emph{symmetric}; their denotations will contain the superscript
$sym$; their elements are called \emph{symmetric functions}. The
next property of the Hoeffding decomposition is specific for the
symmetric case.
\begin{itemize}
\item[iv)] whenever the function $f$ belongs to $
L_p^{sym}(\mu^d),$ the canonical function $ R_S f $ does not
depend on the choice of $S \in \mathcal S^m_d $ and is symmetric;
thus, in this case there exist operators $ R_m: L_p^{sym}(\mu^d)
\to L_p^{sym}(\mu^m)$ such that for every $S=\{i_1, \dots, i_m\}
\in \mathcal S^m_d$
\begin{equation*}
(R_m f) \circ \pi_S = Q_S f.
\end{equation*}
\end{itemize}
Furthermore, every $f \in L_p^{sym}\,(\mu^d)$ can be represented
in a unique way in the form
\begin{equation}\label{symsplit}
f=\sum _{m=0}^d \sum_{S\in\,\mathcal S^m_d}(R_m f)\circ \pi_S.
\end{equation}

\begin{remark}\label{4.3} We illustrate the difference between general and symmetric
  kernels for $d=2$. For a general kernel  $f \in
L_{p\,}(\mu^2)$ we have
\begin{equation*}
f(x_1,x_2)=
f_{\emptyset}+f_{\{1\}}(x_1)+f_{\{2\}}(x_2)+f_{\{1,\,2\}}(x_1,x_2),
\end{equation*}
where $$f_{\emptyset}=\int_{X^2}f(z_1,z_2)\mu(dz_1)\mu(dz_2),$$
$$f_{\{1\}}(x_1)=\int_X f(x_1,z_2)\mu(dz_2) -f_{\emptyset}, \,\,\,f_{\{2\}}(x_2) =\int_X
 f(z_1,x_2)\mu(dz_1)-f_{\emptyset},$$ $$f_{\{1,\,2\}}(x_1,x_2)= f(x_1,x_2) -f_{\{1\}}(x_1) - f_{\{2\}}(x_2)
- f_{\emptyset}.$$ Notice, in order to illustrate the
  notion of canonical kernels, that we have for almost every $x_1, x_2 \in X$,
$$\int_{X} f_{\{1\}}(z) \mu(dz)=0,\int_{X} f_{\{2\}}(z) \mu(dz)=0,  $$
$$ \int_{X} f_{\{1,\,2\}}(z_1,x_2)\mu(d z_1)= \int_{X} f_{\{1,\,2\}}(x_1,z_2)\mu(d
z_2)=0.$$
 For a kernel $f \in L_p^{sym}(\,\mu^2)$ the above
relations reduce to $$f(x_1,x_2)= f_0
+f_1(x_1)+f_1(x_2)+f_2(x_1,x_2),$$ where
$$f_0=\int_{X^2}f(z_1,z_2)\mu(dz_1)\mu(dz_2),$$
$$f_1(x)=\int_X
f(x,z)\mu(dz) -f_0 \Biggl(=\int_X f(z,x)\mu(dz) - f_0\Biggr),$$
 $$f_2(x_1,x_2)=f(x_1,x_2)-f_1(x_1) -f_1(x_2)- f_0.$$
Here
$\int_{X} f_1(z) \mu(dz)=0 $,
 $f_2 \in L_p^{sym}(\,\mu^2)$ and  for almost every $x \in X$ we
have
$$ \int_{X} f_2(z,x)\mu(d z)\Biggl(=\int_{X} f_2(x,z)\mu(d z)\Biggr) =0.$$
\end{remark}

\section{Filtrations. Exactness and Kolmogorov property} \label {5}
In the remaining part  of the paper we deal with distributional
convergence of von Mises statistics for a measure preserving
transformation. Our tool here is a kind of martingale
approximation.
For $d=1$ this approximation goes back to
\cite{Go1969},\,\cite{GoLi1978} and \cite{Ma1978} (in the latter
paper only Harris recurrent Markov chains were considered) and was
developed for higher dimensional random arrays in \cite{Go2009}.\\
The additional structure needed is a filtration compatible with
the dynamics defined by a measure preserving transformation. From
now on we restrict ourselves to a class of measure preserving
transformations of probability spaces, which are \emph{exact}
\cite{Ro1961}.
 Let $T$ be a measure preserving transformation of a
 probability space $(X, \mathcal F, \mu).$ The transformation $T$ defines
 a decreasing filtration $ (T^{-k}
\mathcal F)_{k \ge 0}.$ Exactness of $T$ means that $ \bigcap_{k
\ge 0}T^{-k}\mathcal F = \mathcal N,$ where $\mathcal N$ is the
trivial $\s-$field of the space $(X,\mathcal F, \mu).$ As can
easily be seen, every exact transformation is ergodic. The
standard assumption of the ergodic theory is that $(X,\mathcal F,
\mu)$ is a Lebesgue space in the sense of Rokhlin. Under this
assumption it can be shown that, except for the case of the one
point measure space, the Lebesgue space with an exact
transformation is an atomless measure space, hence, is isomorphic
to the unit interval with the Lebesgue measure. As before,  by $ V^* $  we
denote the adjoint (for $p >1$) and the preadjoint (for $p=1$)  of the operator $V.$ As the operator
$V$ acts as an isometry in all $L_p$ spaces, preserves
constants and positivity, the operator $V^*$ also acts on all
these spaces as a contraction which preserves constants and
positivity. The operator $V^*$ is a particular case of a {Markov
transition operator}.

For every $k \ge 0$ we have the relations  $V^{*k} V^k=I$ and $V^k
V^{*k}=E_k,$ where $I$ is the identity operator and $E_k=
E^{T^{-k}\mathcal F},$ the corresponding conditional expectation.
Let $E$ denote the expectation operator. We can easily conclude
(for example, from known facts about the convergence of reversed
martingales) that the exactness of $T$ is equivalent to the the
strong convergence $V^{*n} \underset{n \to \infty}{\to} E$ in
every space $L_p(\,\mu)$ with $1 \le p < \infty. $ In the sequel
the strong convergence of the series
\begin{equation} \label{series}
 \sum_{k \ge 0} V^{*k}f
\end{equation}
  and other similar conditions will be
imposed on $f.$ Set
$$L_p^0(\mu)=\{f \in L_p(\mu), Ef=0\}.$$
Assuming $T$ is exact, for every $1 \le p < \infty$ the series
\eqref{series}
 converges in the norm of $L_p(\mu)$ if and only if $f$
can be represented in the form $f =(I-V^*)g$ with some $g \in
L_p(\mu)$ (such $g$ is unique up to an additive constant which can
be fixed by the condition $g \in L_p^0(\mu)$). Observe that, in
view of exactness, such $f$'s form a dense subspace in $L_p^0(\mu)$.
\begin{remark} \label{exactness}
In the rest of the paper we will mainly restrict ourselves to
exact transformations. This is just done to simplify the
statements of the results and make the notation more convenient.
We could easily extend these results to ergodic transformations
$T$ and to kernels $f \in L_{p}(\mu^d)$ satisfying the additional
condition $E(f\,|\,\mathcal F_1\otimes \cdots \otimes T_l^{-n}
\mathcal F_l \otimes \cdots \otimes \mathcal F_d)\underset{n \to
\infty}{\to} \check{E}^lf$, $l=1, \dots,d$. Here $T_l$ is the copy
of $T$ acting on the $l$-th coordinate in $X^d$, $\check{E}^l$ was
defined in Subsection \ref{4.1}.
\end{remark}
\begin{remark} \label{adapt}
The results of the next sections are primarily concerned with
exact (hence, non-invertible) transformations; however, they can
be converted into some results on invertible transformations
furnished with an additional structure. Indeed, assume that an
invertible measure preserving $T$ acts on $(X,\mathcal F, \mu)$
and we are given a $\s$-field $ \mathcal F_0 \subset \mathcal F$
such that $T^{-1}\mathcal F_0 \supseteq \mathcal F_0.$ Then a
theory, totally parallel to that we develop in the following
sections for the exact case, applies to kernels measurable with
respect to $\mathcal F_0^{\otimes d}.$ The restriction of $T^{-1}$
to $\mathcal F_0$ corresponds to a non-invertible transformation.
We leave details of this correspondence to the reader; it will be
used when considering applications in Section 9. Just notice that
the counterpart of exactness for an invertible $T$ is the property
$\bigcap_{k \ge 0}T^k\mathcal F_0=\mathcal N$.
 If, moreover,
$\bigvee_{k \ge 0} T^{-k} \mathcal F_0=\mathcal F$, the
transformation $T$ is called \emph{Kolmogorov}. Similarly to the
exactness property in Remark \ref{exactness}, the Kolmogorov
property can be relaxed to the requirement that $T$ is ergodic and
$f$ satisfies an analogue of the additional condition there.
\end{remark}

\section{Growth rates for multiparameter sums} \label{6}

It follows from Lemma \ref{embedding} for $p \in [1, \infty)$ that
the space $ L_{p,\, \pi}^{sym}(\,\mu^m)$ can be identified, using
the injective map $J_m$, with a (non-closed) dense subspace of
$L_{p\,}^{sym}(\,\mu^m)$. As we warned the reader above, the
symbol $J_m$ will be omitted and the relation $
L_{p,\,\pi}^{sym}(\,\mu^m)$ $\subset L_{p\,}^{sym}(\,\mu^m)$ will
be assumed instead of $J_m( L_{p,\, \pi}^{sym}(\,\mu^m))$ $\subset
L_{p\,} ^{sym}(\,\mu^m).$ In particular, it makes sense to speak
of canonical elements of $ L_{p,\,\pi\,}^{sym}(\,\mu^m).$

A noninvertible measure preserving transformation $T$ of a
probability space $(X, \mathcal F, \mu)$ has a natural decreasing
filtration given by $(T^{-n}\mathcal F)_{n  \ge 0}.$ We shall use
the following consequence of the Burkholder inequality.
\begin{lemma} \label{martbound1} For every $p \in [2, \infty)$
there exists a constant $C(p)$ such that for every stationary
sequence $(\xi_n)_{n \in \Z}$ of martingale differences in
$L_p(\mu)$ we have
\begin{equation*}
\biggl|\sum_{k=1}^{n-1}\xi_k\biggr|_p \le C(p)\sqrt{n}\,
\bigl|\xi_0\bigr|_p\,.
\end{equation*}
\end{lemma}

\begin{proof} Let
 $p \in [2, \infty).$
Using the Burkholder inequality (Theorem 9 in \cite{Burk1966}) for
the original sequence and then applying the triangle inequality
for the space $L_{p/2}$ to the sequence $(\xi_n^2)_{n \in \Z}$ ,
we obtain

\begin{equation*}
\frac{1}{\sqrt{n}} \biggl|\sum_{k=1}^{n-1}\xi_k \biggr|_p \le C(p)
\biggl|\biggl(\frac{1}{n}\sum_{k=1}^{n-1}\xi_k^2\biggr)^{1/2}\biggr|_p
\le C(p) |\,\xi_0^2|^{1/2}_{p/2}= C(p)|\,\xi_0|_p\,.
\end{equation*}
\end{proof}
For every $m,$ $0 \le m \le d,$ let  $\mathcal S_m$ ($\mathcal
S_m^s$) be the set of all subsets (res\-pec\-ti\-vely, of all
subsets of cardinality $ s \in \{0,...,m\}\,$) of the set $\{1,
\dots, m\}.$ For every $\mathit{S} \in \mathcal S_m$ define a
subsemigroup $\Z^{m,\mathit{S}}_+ \subseteq \Z^m_+$ by
$$\Z^{m,\mathit{S}}_+ =\{(n_1,\dots,n_m) \in \Z^m_+:n_k = 0\, \,
\text{for all}\,  k \notin \mathit{S}\}.$$ In this section we
write $\mathbf k$ and $\mathbf n$ for $ (k_1,\dots,k_m)$ and
$(n_1,\dots,n_m)$, respectively; the notation $\mathbf k
\boldsymbol{<} \mathbf k'$ ($\mathbf k \boldsymbol{\le} \mathbf
k'$) means that $k_1 < k'_1, \dots, k_m < k_m'$ (respectively,
$k_1 \le k'_1, \dots, k_m \le k_m'$).
\begin{lemma}\label{martbound2}
Let  $m \in \{1,...,d\}$ and let ${\mathbf e}_1, \dots, {\mathbf e}_m $ denote the
standard basis of $\Z_+^m.$ Then, for every real $p \in [2,
\infty)$ and every integer $s \in \{1,...,m\}$,  there exists a
constant $C(p,s)>0$\! with the following property: For every $S
\in \mathcal S_m^s$ and $f \in  L_{\,p,\,\pi}(\,\mu^{m}),$
satisfying
\begin{equation} \label{potential}
 V^{* {\mathbf e}_l}f=0,\, l \in S,
\end{equation}
 the relation
\begin{equation}\label{tensorbound1}
\biggl|\underset{\substack{\mathbf  k \,\in \,\Z^{m,\,\mathit{S}}_+\\
0\,\le\, k_l \,\le\, n_l-1,\,\,  l \,\in\, S}}{\sum} V^{\mathbf
k}f \biggr|_{p,\, m,\,\pi} \le C(p,s) \biggl(\,\prod_{l\, \in \,
S}\sqrt{n_l}\biggr)\, |\,f|_{\,p,\, m,\,\pi}
\end{equation}
holds for every family $(n_l)_{l \in S}$ of natural numbers.
 Moreover, if $p \ge m$ and $r=p/m,$ we  also
obtain
\begin{equation*}
\biggl|\underset{\substack{\mathbf k \,\in\, \Z^{m,\,\mathit{S}}_+\\
0\,\le\, k_l\, \le\, n_l-1,\,\,  l \, \in \, S}}{\sum} D_m
V^{\mathbf k}f \biggr|_r\le C(p,s) \biggl(\,\prod_{l\, \in \,
S}\sqrt{n_l}\biggr)\, |\,f|_{\,p,\,m,\,\pi}
\end{equation*}
for every $(n_l)_{l \,\in \, S}$.
\end{lemma}

\begin{proof} Let $s$ and $S$ be as in the statement
of the lemma. Since the norm of the map $D_m\!\!:
L_{\,p,\,\pi}(\mu^m)\to L_r(\mu)$ is $1$, it suffices to prove
\eqref{tensorbound1}. Let $\mathbf 0_m$ denote the neutral element
of $\Z^m_+$. Set
$$M^S_{\,p,\,\mathbf 0_m,\, \pi}= \{f \in  L_{p,\,\pi}(\,\mu^m):
V^{*{\mathbf e}_l}f=0\, \, \text{for every}\, \, l \in S\}.$$ Observe that
the subspace $M^S_{\,p,\,\mathbf 0_m, \,\pi}\subset
L_{p,\,\pi}(\,\mu^m)$ itself can be  represented as the projective
tensor product of $s$ copies of the subspace
$M_{\,p,\,0}\overset{\mathrm{def}}{=}\{f \in L_p(\,\mu): Vf=0\}$
and $m\!-\!s$ copies of the space $L_p(\,\mu).$ Notice that the
relations \eqref{potential} are equivalent to the following
description of the corresponding subspace in terms of projections:
\begin{equation*}
(I- V^{{\mathbf e}_l} V^{* {\mathbf e}_l})f=f \,\, \text{for every}\, l \in S.
\end{equation*}
 The subspace $M^S_{\,p,\,\mathbf 0_m,\, \pi}$ can
also be
described as the range of the projection
\begin{equation*}
\prod_{l \in\, S}(I- V^{{\mathbf e}_l} V^{* {\mathbf e}_l}).
\end{equation*}
 We need now the following consequence of
Proposition 2.4 in \cite{Ry2002}. In general, for some Banach
spaces $A_l$ and their closed subspaces $B_l\subset A_l,$
$l=1,...,m$, we only have a canonical linear map $ i: B_1
\hat{\otimes}_{\pi}\cdots \hat{\otimes}_{\pi}B_m  \to A_1
\hat{\otimes}_{\pi}\cdots \hat{\otimes}_{\pi}A_m$ of norm $1$.
However, if every $B_l$ is a {complemented} subspace in the
corresponding $A_l$ (that is the range of a bounded projection
$\varphi_l: A_l \to B_l$) then this map is a topological linear
isomorphism onto its range (the latter is closed in $A_1
\hat{\otimes}_{\pi}\cdots \hat{\otimes}_{\pi}A_m $). Moreover, if
every $\varphi_l$ is a projection of norm $1$ then this map is an isometry.\\
Thus, {if  bounded projections} $(\varphi_l)_{l=1, \dots, m}$
exist, we can consider \linebreak $B_1 \hat{\otimes}_{\pi}\cdots
\hat{\otimes}_{\pi}B_m $ as a closed subspace of $A_1
\hat{\otimes}_{\pi}\cdots \hat{\otimes}_{\pi}A_m$, the map
$\varphi_1{\otimes}_{\pi}$ $\cdots {\otimes}_{\pi}\varphi_m $
being a bounded projection of $A_1 \hat{\otimes}_{\pi}\cdots
\hat{\otimes}_{\pi}A_m $ onto its subspace $B_1
\hat{\otimes}_{\pi}\cdots \hat{\otimes}_{\pi}B_m. $ The latter
subspace can be described by
\[B_1\hat{\otimes}_{\pi}\cdots
\hat{\otimes}_{\pi}B_m = \bigl\{f \in A_1
\hat{\otimes}_{\pi}\cdots \hat{\otimes}_{\pi}A_m: (\varphi_1
{\otimes}_{\pi}\cdots {\otimes}_{\pi}\varphi_m )f=f\} \]
or,
equivalently, by
\begin{eqnarray*}
&B_1\hat{\otimes}_{\pi}\cdots \hat{\otimes}_{\pi}B_m = \bigl\{f
\in A_1 \hat{\otimes}_{\pi}\cdots \hat{\otimes}_{\pi}A_m: \\
&\bigl((I-\varphi_1) {\otimes}_{\pi}\cdots {\otimes}_{\pi}I
\bigr)f=0; \dots; \, \bigl(I {\otimes}_{\pi}I{\otimes}_{\pi}
\cdots {\otimes}_{\pi}(I- \varphi_m)\bigr) f=0\bigr\}.
\end{eqnarray*}
Moreover, the projective tensor norm on the space $B_1
\hat{\otimes}_{\pi}\cdots \hat{\otimes}_{\pi}B_m $ and the norm
induced by its embedding
 into $ A_1 \hat{\otimes}_{\pi}\cdots
\hat{\otimes}_{\pi}A_m $ are equivalent.

 We will apply
this assertion to the case when $A_l=L_p(\mu)$ for every $l \in
\{1, \dots, m\},$ $B_l= M_{p\,,0},$ $\varphi_l=I-VV^* $ for $l \in
S,$ and $B_l= L_{p\,}(\mu),$ $\varphi_l=I$ for $l \notin S.$ Since
$VV^*$ is a conditional expectation, it is clear that $\varphi_l$
is bounded for every $l$ (in fact its norm does not exceed
$2^{1-(2/p)}$). With this notation we have that
$M^S_{\,p,\,\mathbf 0_m,\, \pi}$ and $B_1\hat{\otimes}_{\pi}\cdots
\hat{\otimes}_{\pi}B_m$ are isomorphic as topological vector
spaces. Observe that we have here a vector space which is equipped
with  two possibly different norms: the norm inherited from $
L_{\,p\,,\, \pi}(\mu^m)$ and the projective tensor product norm,
respectively. According to one of the properties of the projective
tensor norm
 (\cite{Ry2002}, Proposition 2.8), for every $f \in M^S_{\,p,
 \,\mathbf 0_m,\, \pi}$ and $\e> 0$ there exists a bounded family
 of functions $f_{\,i,\,l} \in B_l$ $(1\le i< \infty, 1 \le l \le m)$
 such that
\begin{equation*}
f =\sum_{i} f_{i,1}\otimes \cdots \otimes f_{i,\,m} \ \mbox{\rm
and}\ \sum_{i} |f_{\,i,\,1}|_p\cdots |f_{\,i,\,m}|_{\,p} \le
C\,'(p,s)|f|_{\,p,\,m,\,\pi} + \e.
\end{equation*}
The constant $C\,'(p,s)$  appears here because we put into the
right hand side the inherited norm $|f|_{\,p,\,m,\,\pi}$ of $f$
rather than its norm in $B_1 \hat{\otimes}_{\pi}\cdots
\hat{\otimes}_{\pi}B_m $.
For $l=1, \dots,m $ and every $i$ let
$F_{\,i,\,l}=\sum_{0 \le\, k \le\, n_l-1} V^k f_{\,i,\,l}$ if $l
\in S,$ and $F_{\,i,\,l}=f_{\,i,\,l}$ if $l \notin S.$ Then,
applying Lemma \ref{martbound1} to the sums $\sum_{k=0}^{n_l-1}V^k
f_{\,i,\,l}$ for $l \in S$ (in this case the summands form a
stationary sequence of reversed martingale differences), it
follows that
\begin{equation*}
\begin{split}
&\biggl|\underset{\substack{\mathbf k \in \,\Z^{m,\mathit{S}}_+
\\ 0\,\le\, k_l \le\, n_l-1,\, \, l \in\, S}}{\sum} V^{\mathbf k}f
\biggr|_{\,\!p,\,m,\,\pi}\\
\le & \sum_{i}\,\Bigl|\!\!\!\!\!\!\!\!\underset{\substack{\mathbf
k \in \,\Z^{m,\mathit{S}}_+ \\0\,\le\, k_l \,\le\, n_l-1,\,\, l
\in \, S}}{\sum} \!\! \!\!\! V^{\mathbf k}\bigl(f_{\,i,\,1}\otimes
\cdots \otimes f_{\,i,\,m}\bigr) \Bigr|_{\,\!p,\,m,\,\pi}\!=\!
\sum_{i}\bigl|F_{\,i,\,1}\otimes \cdots \otimes
F_{\,i,\,m}\bigr|_{\,\!p,\,m,\,\pi}\\
= & \sum_{i}\prod_{l \in  \{1, \dots,m\}}|F_{\,i,\,l}|_p  =
\sum_{i} \prod_{l \in S}\Bigl|\sum_{k=0}^{n_l-1}V^k
f_{\,i,\,l}\Bigr|_p
\prod_{l\notin S}\bigl| f_{\,i,\,l}\bigr|_p  \\
 \le\, & C^s(p)\biggl(\prod_{l \in
S}\sqrt{n_l}\biggr)\sum_{i}\prod_{l \in \,
\{1,\dots,\,m\}}|f_{\,i,\,1}|_{\,\!p}
\cdots |f_{\,i,\,m}|_{\,p}\\
&\le C^s(p)\Bigl(\prod_{l \in
S}\sqrt{n_l}\Bigr)\biggl(C'(p,s)|f|_{\,p_,\,m,\,\pi} + \e\biggr).\\
\end{split}
\end{equation*}
Thus inequality \eqref{tensorbound1} follows with $C(p,s)=
C^s(p)\, C'(p,s).$
\end{proof}
\begin{remark} \label{canonicity} Every $f$ satisfying
the assumptions of the above lemma is $S-$canonical in the
following sense: since every operator $ V^{* {\mathbf e}_l}$ preserves the
integrals with respect to the $l$--th variable, it follows from
\eqref{potential} that, under the assumptions of Lemma
\ref{martbound2}, integrating $f$ over the $l$--th variable
returns $0$ whenever $l \in S.$ This implies the assertion.
\end{remark}
  The following lemma provides a
condition under which the martingale-coboundary decomposition is
valid.
\begin{lemma} \label{martcobdec}
Let  $p \in [1, \infty]$ and  $f \in   L_{p\,,\pi}(\,\mu^m)$ be a
canonical kernel such that the series in the right hand side of
\begin{equation} \label{potential2}
g =\sum_{\mathbf 0 \boldsymbol{\le}\mathbf k \boldsymbol{<}
\boldsymbol{\infty}} V^{* \mathbf k}f \, \,
 \biggl(\,\overset{\mathrm{def}}{=}\,
 \underset{\substack{\,n_1 \,\to\, \infty
\\\, \dots \\\,\, n_m\, \to \, \infty}}{\lim}\,\,\,
\sum_{\mathbf 0 \boldsymbol{\le}\, \mathbf k\boldsymbol{<}\mathbf
n} V^{*\mathbf k}f \biggr)
\end{equation}
converges in $ L_{p\,,\, \pi}(\,\mu^m).$ Then $f$ can be
represented in the form
\begin{equation}\label{representation}
 f=\sum_{S\, \in\, \mathcal S_m} A^Sf,
 \end{equation}
 where for every $S \in \mathcal S_m$
\begin{equation} \label{summand_1}
 A^Sf =  \bigl(\prod_{l\, \notin \,
S}(I-V^{{\mathbf e}_l}V^{*{\mathbf e}_l})\prod_{l\, \in\, S}(V^{{\mathbf e}_l}-I)\bigr)h^S
\end{equation}
and the function $h^S \in  L_{\,p\,,\, \pi}(\,\mu^m)$  is
  defined by the equation
\begin{equation}\label{summand_2}
h^S=\bigl(\prod_{l\, \in \, S}V^{*{\mathbf e}_l}\bigr) g. \,
\end{equation}
The functions $g$ and $(h^S)_{S \in\, \mathcal S_m}$ are
canonical; the summands of the form \eqref{summand_1} in
\eqref{representation} are uniquely determined.
\end{lemma}

\begin{proof} The results and the proofs in \cite{Go2009},
developed originally for the $L_p$--spaces, apply to the
$L_{p\,,\,\pi}$-spaces without any changes. The requirement of
\emph{complete commutativity} imposed in \cite{Go2009} on the
multiparameter dynamical system and the invariant measure is
obviously fulfilled for a direct product with a coordinatewise
action which we deal with in the present paper. Hence, by
Proposition 3 in \cite{Go2009}, the convergence of the series
\eqref{potential2} implies that the \emph{Poisson equation}
(see \cite{Go2009}) is solvable for $f$; therefore, we may apply
Proposition 1 in \cite{Go2009} to $f$. Then we obtain the
representation \eqref{representation} with $A^Sf$ defined by
formulas \eqref{summand_1}, \eqref{summand_2} and the assertion on
the uniqueness of the summands of the form \eqref{summand_1}.
Notice that the operator $V^*$ preserves integrals of functions
with respect to $\mu$; as a consequence, every $V^{*\mathbf n}$
maps canonical functions to canonical ones. Being according to
\eqref{potential2} a limit of canonical functions, $g$ is
canonical. In view of \eqref{summand_2}, all $h^S$ are canonical,
too.
\end{proof}
\begin{proposition} \label{canonkernelbound1} Let
 $0 \le s \le m$ and $f $ be a kernel satisfying the assumptions
of Lemma \ref{martcobdec} for some $p \in [2, \infty).$ Let $A^Sf$
be defined by formulas \eqref{summand_1} and \eqref{summand_2}.
Then there exists a constant $C_{\,p\,,\,m,\,s} >0$\! such that
for every $S \in \mathcal S_m^s$ and every $n_1, \dots, n_m$
\begin{equation}
\label{tensorbound3}\bigl|
\sum_{\mathbf 0 \boldsymbol{\le}\mathbf k \boldsymbol{<} \mathbf
n} V^{\mathbf k}A^Sf \bigr|_{\,p,\,m,\,\pi} \le
C_{\,p,\,m,\,s}\biggl(\prod_{l \notin S}\sqrt{n_l}\biggr)\,
|\,g|_{\,p,\,m,\,\pi},
\end{equation}
where $g$ is defined in \eqref{potential2}. Moreover, for $p \ge m$
\begin{equation}\label{diagbound3}
\bigl|
\sum_{\mathbf 0 \boldsymbol{\le}\mathbf k \boldsymbol{<}\mathbf n}
D_m V^{\mathbf k} A^Sf \bigr|_{\,r}\le
C_{\,p,\,m,\,s}\biggl(\prod_{\,l \,\notin \,S}\sqrt{n_l}\biggr)\,
|\,g|_{\,p,\,m,\,\pi}
\end{equation}
holds with $r=p/m$.
\end{proposition}

\begin{proof}
Setting $\ovS = \{1, \dots,m\}\! \setminus \! S ,$ we have
\begin{equation} \label{coboundary1}
\begin{split}
&\underset{\mathbf 0 \boldsymbol{\le} \mathbf
k\boldsymbol{<}\mathbf n}
{\sum}V^{\, \mathbf k\,} A^S f
\\= \underset
{\substack{\mathbf k\, \in \,\Z^{m,\,\overline{S}}_+\\
0\, \le\, k_t \,\le \, n_t-1, \, t \in
\overline{S}}}{\sum}&V^{\,\mathbf k} \prod_{r \notin
S}(I-V^{{\,\mathbf e}_r}V^{\,*{\mathbf e}_r})
\underset{\substack{\,\mathbf l \, \in \,\Z^{m,\,S}_+\\
0 \,\le\, l_u\, \le\, n_u-1, \, u \in\, S}}{\sum} V^{\,\mathbf
l}\prod_{u\, \in \,
S}(V^{{\,\mathbf e}_u}-I)\, h^{\mathit{S}} \\
= \underset{\substack{\mathbf k \,\in\, \Z^{m,\,\overline{S}}_+\\
0\, \le\, k_t\, \le\, n_t-1, \, t \in
\overline{S}}}{\sum}&V^{\,\mathbf k} \prod_{r\, \notin \,
S}(I-V^{{\,\mathbf e}_r}V^{\,*{\mathbf e}_r})\prod_{u \in
S}(V^{\,n_u{\mathbf e}_u}-I)\,h^{\mathit{S}}.\end{split}
\end{equation}
 Since for $l \notin S$
$$ V^{\,* {\mathbf e}_l}\prod_{r \notin S}(I-V{^{\,\mathbf e_r}}V^{\,*{\mathbf e_r}})
\prod_{u\, \in\,S}(V^{\,\,n_u{\mathbf e}_u}-I)\,h^{\mathit{S}}=0$$
and
$$\Bigl | \prod_{\,r \,\notin \,S}(I-V^{{\mathbf e}_r}V^{\,*{\mathbf e}_r})\prod_{\,u \,\in
S}(V^{\,n_u{\mathbf e}_u}-I)\,h^{\mathit{S}}
\Bigr|_{\,p\,,\,m,\,\pi} \le 2^m \bigl|g\bigr|_{\,p,\,m,\,\pi},
$$
the proposition follows  with $C_{p,\,m,\,s}=2^{m}C(p,s)$ from
Lemma \ref{martbound2}  and formula \eqref{coboundary1}.
\end{proof}

\begin{proposition} \label{canonkernelbound} Let  $p \ge 2$ and  $f \in
L_{\,p\,,\,\pi}(\,\mu^m)$  be a canonical kernel such that the
series on the right hand side of
\begin{equation*}
g =\underset{\mathbf 0 \boldsymbol{\le}\mathbf k \boldsymbol{<}
\boldsymbol{\infty}}{\sum} V^{* \mathbf k}f \, \,
\end{equation*}
converges in $ L_{p,\,\pi}(\mu^m).$ Then for every $n_1, \dots,
n_m$ the following inequality holds
\begin{equation}
\label{tensorbound5}\bigl|\underset{\mathbf 0 \boldsymbol{\le}
\mathbf k\boldsymbol{<}\mathbf n}{\sum} V^{\mathbf k}f
\bigr|_{\,p\,,\,m,\,\pi} \le C_{p\,,\,m}\sqrt{n_1 \cdots n_m}\,
|\,g|_{\,p\,,\,m,\,\pi},
\end{equation}
where $C_{m,\,p}$ is a constant depending only on $m$ and $p$. If,
in addition, $p\in [m, \infty)$ then, with $r =p/m$, we also have
that
\begin{equation*}
\bigl|\underset{\mathbf 0 \boldsymbol{\le} \mathbf
k\boldsymbol{<}\mathbf n}{\sum}  D_m V^{\mathbf k}f
\bigr|_{\,r}\le C_{p\,,\,m}\sqrt{n_1 \cdots n_m}\,
|\,g|_{\,p\,,\,m,\,\pi}.
\end{equation*}
\end{proposition}

\begin{proof} Again, since the norm of the operator $D_m:
L_{p\,,\,\pi}(\mu^m) \to L_r(\mu)$ is $1$, we only need to prove
\eqref{tensorbound5}. As $n_1 \ge 1, \dots, n_m \ge 1,$ we have
for every $S \in \mathcal S_m$ $ \prod_{\,l\, \in \,S}
\frac{1}{n_l} \le 1$. Using this relation along with
\eqref{representation} and \eqref{tensorbound3} we obtain
\eqref{tensorbound5} with $C_{p\,,\,m}=\sum_{s=0}^m {m\choose
s}C_{\,p,\,m,\,s}$.
\end{proof}
The following sufficient condition for convergence of the series
in \eqref{potential2} will be used in Section 9 when considering
applications. Expansion of
a kernel into an absolutely
convergent series whose summands are products of functions in
separate variables is
natural in the context of the limit
theory of $U$- and $V$-statistics (see, for example,
\cite{BoVo2008}). Projective tensor products
call for using such
series
to representing
arbitrary elements (see Proposition
2.8 in \cite{Ry2002}). Neither
uniqueness of the representation, nor
 linear independence of the 'basis' is assumed. Notice that
we used such a decomposition in Corollary \ref{sllnser}.

\begin{proposition} \label{suffic} Let, for some $p \in [1,\infty]$,
$(e_k)_{k=0}^{\infty}$ be a sequence of functions such that
$e_0\equiv 1$ and  for every $k \ge 1$ $e_k \in L_p\,(\mu)$ with
$\int_X e_k(x) \mu(dx)=0$. Assume that for every $k \ge 1$
$$C_{p\,,\,k}\eqdef \sum_{n \, \ge\, 0}|\,V^{*n} e_k\,|_p <\infty. $$
Suppose that $f \in L_{\,p}(\,\mu^m)$
admits a representation
\begin{equation} \label{multiple}
f(x_1,\dots,x_m)=
\sum_{\mathbf 0 \boldsymbol{<} \mathbf k
\boldsymbol{<}\boldsymbol{\infty}}
\l_{\bold k}(f)\, e_{k_1}(x_1)\cdots e_{k_m}(x_m)
\end{equation}
where $(\l_{\mathbf k}(f))_{\mathbf 0
\boldsymbol{<}\mathbf k \boldsymbol{<}\boldsymbol{\infty}}$
is a family of constants satisfying
\begin{equation} \label{up}
 C_{p\,}(f)\eqdef\sum_{\mathbf 0
\boldsymbol{<}\mathbf k \boldsymbol{<}\boldsymbol{\infty}}
|\,\l_{\,\mathbf k}\,(f)|\, C_{p\,,\,k_1}\cdots\, C_{p\,,\,k_m} <
\infty.
\end{equation}
 Then $f$ is a canonical kernel of degree $m$,
$f \in L_{p,\,\pi}(\,\mu^m)$, the series in \eqref{potential2}
converges in $L_{p\,,\,\pi}(\,\mu^m)$ and its sum $g$ satisfies
the
inequality
\begin{equation} \label{pronorm}
|\,g|_{\,\!p\,,\,m,\,\pi}\le C_{p\,}(f).
\end{equation}
\end{proposition}
\begin{proof}
For every $k\ge 1$ $C_{p\,,\,k} \ge |\,e_k|_p\,$. Hence,
$$\sum_{\mathbf 0
\boldsymbol{<}\mathbf k \boldsymbol{<}\boldsymbol{\infty}}
|\,\l_{\,\mathbf k}\,(f)|\,|\,e_{k_1}\,|_{\,p}\cdots
|\,e_{k_m}\,|_{\,p} \le C_{p\,}(f) < \infty.$$ Then, according to
\cite{Ry2002}, $|f|_{\,p\,,\,m,\,\pi} \le C_{\,p}(f) < \infty$ and
$f \in L_{p\,,\,\pi\,}(\,\mu^m)$; $f$ is canonical because so is
every term of the series in \eqref{multiple}. Now we
obtain
\begin{equation*}
\begin{split}
&|\,g|_{\,p,\,m,\,\pi} \le \sum_{\mathbf 0 \boldsymbol{\le}\mathbf
n \boldsymbol{<} \boldsymbol{\infty}}\,\,\, \sum_{\mathbf 0
\boldsymbol{<}\mathbf k \boldsymbol{<}\boldsymbol{\infty}}
|\,\l_{\,\mathbf k}\,(f)|\, |V^{* \mathbf n} (e_{k_1}\cdots
e_{k_m})|_{\!\,p,\,m,\,\pi}\\
&= \sum_{\mathbf 0 \boldsymbol{\le}\mathbf n \boldsymbol{<}
\boldsymbol{\infty}}\,\,\, \sum_{\mathbf 0 \boldsymbol{<}\mathbf k
\boldsymbol{<}\boldsymbol{\infty}} |\,\l_{\,\mathbf k}\,(f)|\,
|V^{* n_1} e_{k_1}|_p\cdots |V^{*n_m }e_{k_m}|_{\,p}\\
&= \sum_{\mathbf 0 \boldsymbol{<}\mathbf k
\boldsymbol{<}\boldsymbol{\infty}}|\,\l_{\,\mathbf k}\,(f)|\,
C_{p,\,k_1}\cdots\, C_{p,\,k_m}= C_p(f)< \infty.
\end{split}
\end{equation*}
\end{proof}

\section{Central Limit Theorems in the non-degenerate case} \label{7}

$N(m, \s^2)$ will denote the Gaussian distribution in $\R$ with
mean value $m \in \R$ and variance $\s^2\ge 0$ including the case
$\s^2=0$ of the Dirac measure at $m \in \R$. We first prove a
central limit theorem together with the convergence of the second
moments.
\begin{theorem} \label{clt2}
Let $f \in  L_2^{sym}(\,\mu^d)$ be a real valued kernel with the
symmetric Hoeffding decomposition
\begin{equation*}
f=\sum _{\,m\,=\,0\,}^d \sum_{\,S\,\in\,\mathcal S^m_d}(R_m
f)\circ \pi_S.
\end{equation*}
Assume that
 for every $m=1, \dots, d$ $R_m f \in
 L_{2m,\,\pi\,}^{sym}(\mu^m)$ and that the series
\begin{equation} \label{thesum2a}
\sum_{\substack{\mathbf k \in \Z^m_+\\\mathbf 0
\boldsymbol{\le}\mathbf k \boldsymbol{<} \boldsymbol{\infty}}}
V^{\,* \mathbf k}\,R_m f \,
\biggl(\,\overset{\mathrm{def}}{=}\,\,\,\underset{\substack{n_1 \to \infty\\ \dots \\
n_m \to \infty}}{\lim} \sum_{\substack{\mathbf k \in
\Z^m_+\\\mathbf 0 \boldsymbol{\le}\mathbf k \boldsymbol{<} \mathbf
n}} V^{\,* \mathbf k} R_m f \biggr)
\end{equation}
 converges in $
L_{2m,\, \pi}(\mu^m).$
 Then the sequence
$$
V_n^{(d)}f = \frac{1}{n^ {d-1/2}} \underset{\substack{\,0\le k_1
\le n-1
\\\, \dots \\\,\, 0 \le k_d \le n-1}}{\sum} D_d V^{(k_1,\, \dots,\,k_d)}(f-R_0 f)$$
converges in distribution to $N(0, {d\,}^2 \sigma^2(f)),$ where
\begin{equation*}
\sigma^2(f)= \biggl|\sum_{k=0}^{\infty} V^{* k} R_1 f \biggr|_2^2
- \biggl|\sum_{k=1}^{\infty} V^{* k} R_1 f\biggr|_2^2\ge 0.
\end{equation*}
The convergence of the second  moments
\begin{equation*}
E(V_n^{(d)}f)^2\underset{n \to \infty}{\to}{{d\,}^2 \sigma^2(f)}
\end{equation*}
holds as well.
\end{theorem}

\begin{remark}
 According to the standard terminology, a kernel $f$ is called
\emph{non-degenerate}  if $R_1f$ does not vanish  identically,
otherwise $f$ is called \emph{degenerate}. In the case of
i.i.d.~variables such non-degeneracy is equivalent to the
non-degeneracy of the limit Gaussian distribution using 
normalization by the constants $n^ {d-1/2}$. However, in the
general stationary dependent case such a \emph{statical
non-degeneracy} may occur together with the degeneracy of the
limit distribution. This phenomenon can be viewed as a
\emph{dynamical degeneracy}.
\end{remark}

\begin{proof}
Decompose $f-R_0 f$ in the following way:
\begin{equation*}
f-R_0 f=\sum _{\,m\,=\,1\,}^d \sum_{S\,\in\,\mathcal S^m_d}(R_m
f)\circ \pi_S=\sum _{m\,=\,1}^d f_m,
\end{equation*}
where
$$f_m=
\sum_{S\,\in\,\mathcal S^m_d}(R_m f)\circ \pi_S,\qquad m=1, \dots,
d.
$$
In order to prove the theorem it suffices to establish that
\begin{trivlist}
\item 1) $ V_n^{(d)}f_1 $ converges in distribution to
$N(0,d^2\sigma^2(f)),$ \item 2) $|V_n^{(d)}f_1|_{\,2}^{\,2}
\underset{n \to\infty}{\to}d^2\sigma^2(f),$
 \item 3)
$|V_n^{(d)}\sum_{m=2}^df_m|_{\,2} \underset{n \to \infty}{\to} 0.$
\end{trivlist}
In view of the equality
\begin{equation*} \label{recip2}
D_d \Bigl(V^{(\,k_1,\,\dots,\,k_d)}\sum_{S\in\mathcal S^m_d}(R_m
f)\circ \pi_S\Bigr) = \!\sum_{S=\{i_1,\,\dots,\,i_m\}\in\,\mathcal
S^m_d}D_m V^{(\,k_{i_1},\dots,\,k_{i_m})}R_m f
\end{equation*} we obtain
\begin{equation} \label{reduct2-2}
\begin{split}
&V_n^{(d)}f_m=V_n^{(d)}\Bigl(\sum_{S\in\mathcal S^m_d}(R_m f)\circ
\pi_S\Bigr)=\\
 &=\frac{1}{n^ {d-1/2}} \underset{\substack{0\le k_1 \le n-1
\\ \dots \\ 0 \le k_d \le n-1}}{\sum}
D_d\Bigl(V^{(\,k_1,\dots,\,k_d)}\sum_{S\,\in\,\mathcal S^m_d}(R_m
f)\circ
 \pi_S\Bigr)\\
&= \frac{1}{n^ {d-1/2}} \underset{\substack{0\,\le\, k_1\, \le\,
n-1
\\ \dots \\ 0\, \le\, k_d \,\le\, n-1}}{\sum}\,
\, \sum_{\,\,\,S=\{\,i_1,\dots,\,i_m\,\}\,\in\,\mathcal
S^m_d}D_mV^{(\,k_{\,i_1},\dots,\,k_{i_m})} R_m f \\
&= \frac{ {d\choose m}}{n^ {m-1/2}}\,\, D_m
\underset{\substack{0\,\le\, k_1 \,\le\, n-1
\\ \dots \\ 0\, \le\, k_m \,\le\, n-1}}{\sum}\,  V^{(\,k_1,\dots,\,k_m)}R_m f
\end{split}
\end{equation}
for every $m=1, \dots, d.$  It follows from \eqref{reduct2-2},
Proposition \ref{canonkernelbound} with $p=2m$ and the assumptions
of the theorem that the function $f_m$ satisfies the inequality
$$|V_n^{(d)}f_m|_2 \le C_m   {d\choose m} n^{-(\,m-1)/2}\,
|\,g_m|_{\,2\,m,\,m,\,\pi}$$ where $g_m$ denotes
the sum of the series \eqref{thesum2a}.
This bound for $m \ge 2$ proves 3).\\
Consider now the sums involving $f_1.$ We obtain from
\eqref{reduct2-2} that
\begin{equation} \label{transform}
V_n^{(d)}f_1= d
\frac{1}{\sqrt{n}}\sum_{k=0}^{n-1} V^k R_1 f,
\end{equation}
 where $R_1f$ has the representation
$R_1 f=g_1 - V^{*}g_1$ with $g_1$ denoting the series
\eqref{thesum2a} for $m=1$. This representation can be rewritten
as
\begin{equation} \label{decomposition_2}
R_1 f= (I-V V^*)g_1 + (V-I)V^*g_1.
\end{equation}
 Here the first summand gives, under the action of the operators $(V^k)_{k\, \ge\, 0}$,
 an ergodic stationary sequence of reversed square integrable martingale differences
$(V^k(I-V V^*)g_1)_{k \ge 0}$. By the Billingsley-Ibragimov CLT
\cite{Bi1961,Ib1963}, the variables
$1/\sqrt{n}\,\,\sum_{k=0}^{n-1} V^k (I-V V^*)g_1 $ converge in
distribution, along with the variance, to the required  centered
Gaussian law. The second summand in \eqref{decomposition_2} only
makes a uniformly $L_2$--bounded contributions to each of the sums
$\sum_{k=0}^{n-1} V^k R_1 f.$ Thus, the convergence to the Gaussian distribution in 1) is established. \\
The convergence of the second moments can be concluded as follows.
In the situation of the Billingsley-Ibragimov CLT   we have
$$\biggl| \frac{\sum_{k=0}^{n-1} V^k (I-V V^*)g_1}{\sqrt n}\biggr|_2^2=
|(I-V V^*)g_1|_2^2=|g_1|_2^2-|V^*g_1|_2^2 = \sigma^2(f).$$ This
implies, in view of \eqref{transform}, \eqref{decomposition_2} and
the triangle inequality, that
$$ ||V_n^{(d)}f_1|_2-d\sigma(f)| \le \frac{2d |g_1|_2}{\sqrt{n}},$$
which proves 2) and, together with 3), the convergence of the
second moments.
\end{proof}
Under somewhat weaker assumptions we have the
following central limit theorem with the convergence of the first
absolute moment.

\begin{theorem} \label{clt1}
Let $f \in  L_1^{sym}(\,\mu^d)$ be a real valued kernel with the
symmetric Hoeffding decomposition
\begin{equation*}
f=\sum _{m\,=\,0}^d \sum_{S\in\,\mathcal S^{\,m}_d}(R_m f)\circ
\pi_S.
\end{equation*}
Assume that
\begin{enumerate}
 \item  for every $m=1, \dots, d$ $R_m f \in
 L_{m,\,\pi}^{sym}(\,\mu^m)$ and the series
 \begin{equation} \label{thesum1}
\sum_{\substack{\mathbf k \in \Z^m_+\\\mathbf 0
\boldsymbol{\le}\mathbf k \boldsymbol{<} \boldsymbol{\infty}}}
V^{* \mathbf k}\,R_m f \,
\biggl(\,\overset{\mathrm{def}}{=}\,\,\,\underset{\substack{n_1 \to \infty\\ \dots \\
n_m \to \infty}}{\lim} \sum_{\substack{\mathbf k \in
\Z^m_+\\\mathbf 0 \boldsymbol{\le}\mathbf k \boldsymbol{<} \mathbf
n}} V^{* \mathbf k} R_m f \biggr)
\end{equation}
converges in $ L_{m,\, \pi}(\,\mu^m)$, \item
 $R_1 f$ satisfies the relation
\begin{equation} \label{sum_1}
\bigl|\sum_{ k\,=\,0}^{n-1} V^kR_1 f\bigr|_1
=O(\sqrt{n})\qquad\mbox{as $n \to \infty$}.
\end{equation}

\end{enumerate}
Then there exists $\sigma^2(f)\ge 0$ such that the sequence
$$
V_n^{(d)}f = \frac{1}{n^ {d-1/2}} \underset{\substack{\,0\,\le\,
k_1\, \le\, n-1
\\\, \dots \\\,\, 0\, \le\, k_d \,\le \, n-1}}{\sum}
D_dV^{(\, k_1,\dots,\, k_d\, )}(f-R_0 f), \,\,\, n \ge 1,$$
converges in distribution to $N(0, {d\,}^2 \sigma^2(f))$ as $n \to
\infty$. The convergence of the first absolute moments
\begin{equation} \label{firstmoment}
E|\,V_n^{(d)}f|\underset{n \to
\infty}{\to}{d\sqrt{\frac{2}{\pi}}}\,\sigma(f)
\end{equation}
holds as well.
\end{theorem}

\begin{proof}
The proof is parallel to that of Theorem \ref{clt2}, so we will
concentrate on the essential changes in the proof. Consider the
Hoeffding decomposition of $f-R_0 f$
\begin{equation*}
f-R_0 f=\sum _{m=1}^d \sum_{S\in\mathcal S^m_d}(R_m f)\circ
\pi_S=\sum _{m=1}^d f_m
\end{equation*}
with
$$f_m=
\sum_{S\in\mathcal S^m_d}(R_m f)\circ \pi_S,\, m=1, \dots, d. $$
In order to prove the theorem it suffices to establish that
\begin{trivlist}
\item 1) for some $\sigma(f)\ge 0$, $ V_n^{(d)}f_1 $ converges in
distribution to $N(0,d^2\sigma^2(f)),$ \item 2) $|V_n^{(d)}f_1|_1
\underset{n \to\infty}{\to}{d\sqrt{\frac{2}{\pi}}}\,\sigma(f),$
 \item 3)
$|V_n^{(d)}\sum_{m=2}^df_m|_1 \underset{n \to \infty}{\to} 0.$
\end{trivlist}
Analogously to the proof of Theorem \ref{clt2}, the functions
$f_m,$ $ 1 \le m \le d,$ can be shown to satisfy the inequality
$$|V_n^{(d)}f_m|_1 \le C_m   {d\choose m} n^{-(m-1)/2}\,
|\,g_m|_{\,m,\,m,\,\pi},$$  where $g_m$ denotes the sum of the series
\eqref{thesum1}. For $m \ge 2$ the latter bound implies the
convergence  in $L_1(\mu)$ to zero, proving 3). Taking $m=1$, we
obtain
$$ V_n^{(d)}f_1= d \frac{1}{\sqrt{n}}\sum_{k=0}^{n-1} V^k R_1 f,$$
where $R_1f$ has the representation $R_1 f=g_1 - V^{*}g_1$ with
$g_1 \in L_1(\mu)$ denoting the sum of the series \eqref{thesum1}
with $m=1$. As in the proof of Theorem 2, $R_1 f$ can be
represented in the form
\begin{equation*}
R_1 f= (I-V V^{*})g_1 + (V-I)V^*g_1,
\end{equation*}
 where the first summand defines an ergodic stationary
sequence of reversed martingale difference $(V^k(I-V V^*)g_1)_{k
\ge 0}$, while the second one only contributes  a uniformly
$L_1$-bounded amount to each of the sums $\sum_{k=0}^{n-1} V^k R_1
f$. However, now we only have $(I-V V^*)g_1\in L_1(\mu),$ while we
need $(I-V V^*)g_1\in L_2(\mu)$ to apply the Billingsley-Ibragimov
CLT. The latter can be concluded, as suggested in \cite{Go1973},
from \eqref{sum_1} using another Burkholder inequality  (Theorem 8
in \cite{Burk1966}) and the ergodic theorem (see \cite{Brad1988}
for details). This proves the convergence in distribution. The
convergence of the first moments can be concluded similarly to the
corresponding part in the proof of Theorem \ref{clt2}.
\end{proof}
\begin{remark} \label{alrernat} In the statement of Theorem \ref{clt1}
the requirement \eqref{sum_1} can be substituted by the  relation
\[
\Bigl| \underset{\substack{\,0\,\le\, k_1\, \le\, n-1
\\\, \dots \\\,\, 0\, \le\, k_d \,\le \, n-1}}{\sum}
D_dV^{(\, k_1,\dots,\, k_d\, )}(f-R_0 f)\Bigr|_{\,1} =O(
n^{d-1/2})\qquad\mbox{as $n \to \infty$}.
\]
\end{remark}

\section{A limit theorem for canonical kernels of degree $2$} \label{8}

Apart from non-degenerate kernels of the previous section, a
different type of von Mises statistics emerges from canonical
symmetric kernels of degree $d \ge 2$. Limit distributions of
$V$-statistics defined by such kernels are usually described in
terms of series (or polynomials) in Gaussian variables, or in
terms of multiple stochastic integrals. In the case of
$V$-statistics of dependent variables some descriptions of the
limits in terms of dependent Gaussian variables or non-orthogonal
stochastic integrals are known \cite{BoBy2006-1, BoBy2006-2,
Ea1979}. A rather attractive way is to present the limit
distribution, like in the i.i.d. case, in terms of
\emph{independent} Gaussian variables. This will be done below in
the case $d=2$ and is based on the \emph{diagonalization} of the
symmetric kernel. The point is that the diagonalization here is
applied, instead of the original kernel, to a \emph{martingale
kernel} which emerges as a leading summand in the
\emph{martingale-coboundary representation} of the original
kernel. Notice that the diagonalization of martingale kernels is
also used in \cite{LeNe2011}; in the present work, however,
martingale kernels are considered as a subclass
 to which the study of much more general kernels is reduced.\\
We assume that $f=f_2$ in terms of the Hoeffding decomposition for
symmetric kernels (see Remark \ref{4.3} in Subsection \ref{4.2}).
Let $\theta$ denote the involution in $(X^2\!,\mathcal F^{\otimes
2}\!\!,\mu^2)$ interchanging the multipliers in the Cartesian
product. We consider the spaces $ L_{\,2,\,\pi}(\,\mu^2)$ and
$L_{\,2,\,\pi}^{sym}(\,\mu^2)$ as embedded in $L_2(\,\mu^2).$
\begin{proposition} \label{twoprop}  Let
$f(=f_2) \in L_{2,\,\pi}^{sym}(\,\mu^2)$ be a canonical kernel of
degree $2.$ If the limit
\begin{equation} \label{conv8}g \,\overset{\mathrm{def}}{=}\,
\underset{n_1,\,n_2 \to \infty}{\lim}\underset{\substack{0 \, \le\, i_1 \,\le \, n_1-1\\
0\, \le \,i_2\, \le \, \, n_2-1}}{\sum} V^{* (i_1,\,i_2)}f
\end{equation}
 exists in $ L_{2, \pi}(\mu^2)$,
then $f$ admits a unique representation of the form
\footnote{Upper indices here follow Lemma \ref{martcobdec}.}
\begin{equation}\label{decomptwo}
\begin{split} &f \\
&= g^{\emptyset}\! +\!(V^{(1,0)}\!-\!I)g^{\{1\}}
+(V^{(0,1)}\!-\!I)g^{\{2\}}
+(V^{(1,\,0)}\!-\!I)(V^{(0,1)}\!-\!I)g^{\{1,2\}},
\end{split}
\end{equation}
 where
\begin{equation} \label{condit}
\begin{split}
&E(g^{\emptyset}|T^{-(1,\,0)}\mathcal F^{\,\otimes 2})
=0,\,\,\,\,\,\,\,\,\,\,\,
E(g^{\emptyset}|T^{-(0,\,1)}\mathcal F^{\,\otimes 2})=0,\\
&E(g^{\{1\}}|T^{-(0,\,1)}\mathcal F^{\,\otimes 2})=0,\,\,\,
E(g^{\{2\}}|T^{-(1,\,0)}\mathcal
 F^{\,\otimes 2})=0
\end{split}
\end{equation}
 and $g^{\emptyset}$, $g^{\{1\}}$\!, $g^{\{2\}}$,
 $g^{\{1,2\}}$ are canonical. The functions $g^{\emptyset}$, $g^{\{1\}}$\!, $g^{\{2\}}$,
 $g^{\{1,2\}}$ in \eqref{decomptwo} are uniquely determined by the above properties;
 moreover, $\,g^{\emptyset},g^{\{1,2\}}\!\in\! L_{2,\pi}^{sym}(\mu^2)$,
$g^{\{1\}}\!,g^{\{2\}}\!\in \! L_{2, \pi}(\mu^2),$
$g^{\{1\}}\circ\theta\!=\!g^{\{2\}}$ and $g^{\{2\}}\circ\theta=
g^{\{1\}}$.
\end{proposition}
\begin{proof}
Up to the details related to symmetry the proposition follows from
Example 2.1 in \cite{Go2009}. We propose, however, a partially
independent proof based on the decomposition of $f$ presented by
Lemma \ref{martcobdec} with $m =2.$  Set ${\mathbf e}_1=(1,0),
{\mathbf e}_2=(0,1)$. Vanishing of conditional expectations
follows from the presence of the operators $I-V^{\mathbf
e_l}V^{*\mathbf e_l}$ $(l=1,2)$ in corresponding summands of
\eqref{summand_1} (recall that $V^{\,{\mathbf e}_1}$,
$V^{\,*{\mathbf e}_1}$ commute with $V^{\,\mathbf e_2}$, $V^{\,*
\mathbf e_2}$). Further, the operators $I-V^{\,\mathbf
e_l}V^{*{\,\mathbf e}_l}$ preserve canonicity since so do $I,
V^{{\,\mathbf e}_l}$ and $V^{\,*{\mathbf e}_l}$. Hence, the
functions $g^{\emptyset}$, $g^{\{1\}}$, $g^{\{2\}}$, $g^{\{1,2\}}$
are canonical because so are the functions $h^S$ in Lemma
\ref{martcobdec}; this lemma also implies the uniqueness of the
summands in the representation \eqref{decomptwo}. To establish the
uniqueness claimed in the proposition we need to prove that
canonical solutions to the equations
$$(V^{{\,\mathbf e}_1}\!-\!I)g^{\{1\}}=0,
(V^{\,{\mathbf e}_2}\!-\!I)g^{\{2\}}=0, (V^{\,{\mathbf
e}_1}\!-\!I)(V^{\,{\mathbf e}_2}\!-\!I)g^{\{1,2\}}=0$$ vanish.
Applying $V^{\,*\,{\mathbf e}_l}$ to the first equation, we obtain
$(I-V^{*{\,\mathbf e}_l})g^{\{1\}}=0$ or $g^{\{1\}}=V^{*{\,\mathbf
e}_l}g^{\{1\}}$. Iterating the latter equation gives
$g^{\{1\}}=V^{\,*n{\mathbf e}_1}g^{\{1\}}$ for every $n \ge 1$.
 For a canonical $g^{\{1\}}$ the right hand side of
the last equation tends to $0$ as $n\! \to\! \infty$; hence
$g^{\{1\}}=0$. Other equations can be treated similarly. The
symmetry of $g^{\emptyset}$ and $g^{\{1,2\}}$ follows from the
symmetry of $f$ and the uniqueness. Then we apply $\theta$ to the
decomposition of a symmetric $f$ and use the uniqueness of
summands in the decomposition \eqref{decomptwo} with symmetric
$g^{\emptyset}$ and $g^{\{1,2\}}$. By uniqueness we obtain
$g^{\{1\}}\circ\,\theta\!=\!g^{\{2\}}$ and
$g^{\{2\}}\circ\,\theta= g^{\{1\}}$.
\end{proof}

Assume, in addition, that the kernel $f $ is real-valued. The
function $g^{\emptyset}$ is the real-valued kernel of a symmetric
trace class integral operator in $L_2(\mu).$ Hence, it admits the
eigenfunction decomposition
\begin{equation}\label{kernexpand}
 g^{\emptyset}(x_1,x_2)= \sum_{m=1}^{\infty}\lambda_m\varphi_m(x_1)\varphi_m(x_2)
\end{equation}
 where $(\varphi_m)_{m \ge 1}$ is a normalized orthogonal sequence in
$L_2(\mu)$ and \linebreak $(\lambda_m)_{m \ge 1}$ is a real
sequence (of not necessarily distinct numbers) for which
$\sum_{m=1}^{\infty} |\,\lambda_m| < \infty$. We shall assume that
$\lambda_m \neq 0$ for every $m \ge 1$. Moreover, since we
consider $L_2(\mu)$ over $\C$, we assume that the functions
$(\varphi_m)_{m \ge 1}$ are chosen real-valued (this is always
possible since $g$ is symmetric and real-valued).

\begin{theorem}\label{degenCLT} Let $f$ be a real-valued canonical kernel satisfying
the assumptions of Proposition \ref{twoprop}. Then, as $n \to
\infty,$ the sequence of random variables
$$
 \, \frac{1}{n}\,\sum_{0 \le\, i_1\!, \,i_2 \le\, n-1}  D_2 V^{(i_1,\,i_2)} f,\, \, \,\, \,n \ge 1, $$ converges in
distribution to
$$\xi\,\overset{\mathrm{def}}{=}\,
\sum_{m=1}^{\infty} \lambda_m \eta_m^2,$$ where $
(\eta_m)_{m=1}^{\infty}$ is a sequence of independent standard
Gaussian variables. Moreover,
\begin{equation*}
E \,\Bigl(\, \frac{1}{n}\,\sum_{0 \le\, i_1, \,i_2\, \le \,n-1}
D_2\,V^{(i_1,i_2)} f \Bigr) \underset{n \to \infty}{\to}
\sum_{m=1}^{\infty}\lambda_m .
\end{equation*}
\end{theorem}

\begin{proof}
Setting in \eqref{diagbound3} $m=2$,\,$p=2$,\, $r=1$, we obtain
with $s=1$
\[
\Bigl|\sum_{0 \le\, i_1\!, \,i_2 \le\, n-1}\!\!\!\!\!\!\! D_2
V^{(i_1,\,i_2)} \bigl((V^{(1,\,0)}\!-\!I)g^{\{1\}}
+(V^{(0,1)}\!-\!I)g^{\{2\}}\bigr)\Bigr|_1\! \le \!
2\,C_{\,2,\,2,\,1}\sqrt{n}\, |\,g|_{\,2,\,2,\,\pi}
\]
 and with $s=2$
 \[
 \Bigl|\sum_{0 \le\, i_1, \,i_2 \le\, n-1}\!\!\!\!\!\!\! D_2
V^{(i_1,\,i_2)}\bigl((V^{(1,0)}\!-\!I)(V^{(0,1)}\!-\!I)g^{\{1,\,2\}}\bigr)\Bigr|_1
\le  C_{\,2,\,2,\,2} \, |\,g|_{\,2,\,2,\,\pi}.
\]
These two inequalities and decomposition \eqref{decomptwo} imply
that
\begin{equation*}\label{}
\Bigl|\frac{1}{n}\, \sum_{0 \le\, i_1\!, \,i_2 \le\, n-1} D_2
V^{(i_1,\,i_2)} (f-g^{\emptyset})\Bigr|_1\,\, \underset{n \to
\infty}{\to} 0
\end{equation*}
which reduces the proof to the special case of the
kernel $g^{\emptyset}.$\\
 Let us show next that for every $m \ge 1$
$$E(\varphi_m|T^{-1} \mathcal F)=0.$$
We have $\mu\times\mu\,$-almost surely
\begin{equation*}\label{}
0=E(g^{\emptyset}|\,T^{-(0,1)}\mathcal F^{\otimes 2})(x_1,x_2)=
\sum_{l=1}^{\infty}\lambda_l\varphi_l(x_1)E(\varphi_l|\,T^{-1}\mathcal
F)(x_2)
\end{equation*}
which for every $m$ implies, via multiplying by $ \varphi_m(x_1)$
and integrating  over $x_1$ with respect to $\mu$, that $\lambda_m
E(\varphi_m|\,T^{-1}\mathcal F)(x_2)=0$.  Thus we have
\begin{equation*}\label{}
 E(\varphi_m|\,T^{-1}\mathcal F)=0
\end{equation*}
 $\mu-$almost surely for every $m \ge 1$.\\
Define now a random variable $\xi_N$ and a truncated kernel
$g_N^{\emptyset}$ by setting
$$\xi_N=\sum_{m=1}^{N}\lambda_m
\eta^2_m,\,\,\,\,\,g_N^{\emptyset}(x_1,x_2)=
\sum_{m=1}^{N}\lambda_m\varphi_m(x_1)\varphi_m(x_2).$$  Observe
that for every $N$ the assertions of the theorem on  the
convergence in distribution and the convergence of the first
moments hold for  for $g_N^{\emptyset}$ and $\ \xi_N$, when $f$ is
replaced by $g_N^{\emptyset}$. Indeed, the Billingsley-Ibragimov
theorem applies to reversed $\R^N$-valued martingale differences
(this  is straightforward via the Cramer-Wold device). So, the
random vectors
$$ \Bigl(\frac{1}{\sqrt{n}}\sum_{k=o}^{n-1}\varphi_1\circ T^k, \dots,
\frac{1}{\sqrt{n}}\sum_{k=o}^{n-1}\varphi_N\circ T^k\Bigr)$$
converge in distribution to $(\eta_1,\dots,\eta_N)$ as $n \to
\infty.$ Hence, the random variables
$$\frac{1}{n}\,
 \sum_{0 \le\, i_1\!, \,i_2 \le\, n-1} D_2
V^{(i_1,\,i_2)}g^{(N)}_{\emptyset}=
\sum_{m=1}^{N}\lambda_m\biggl(\frac{1}{\sqrt{n}}\sum_{k=0}^{n-1}\varphi_m
\circ T^{\,k}\biggr)^2 $$ converge in distribution to
$\sum_{m=1}^N \l_m\eta_m^2$ as $n \to \infty.$ The convergence of
the first moments follows here from the convergence of the second
moments in the CLT for martingale differences. Observe now that
\begin{equation*}\label{}
|\,\xi-\xi_N|_1 = \Bigl|\sum_{m=N+1}^{\infty}\lambda_m \,
\eta^2_m\Bigr|_1 \le \sum_{m=N+1}^{\infty}|\,\lambda_m|
\underset{N\to\infty}{\to} 0.
\end{equation*}
Hence, $(\xi_n)_{n\ge 1}$ converges to $\xi $ in distribution
along with the first moment. Combining this with the fact that
\begin{equation*}\label{}
\begin{split} \biggl|\frac{1}{n}\,&  \sum_{0 \le\, i_1\!, \,i_2 \le\, n-1} D_2 V^{(\,i_1,\,i_2\,)}g^{\emptyset}-\frac{1}{n}\, \sum_{0 \le\,
i_1\!, \,i_2 \le\, n-1} D_2
V^{(\,i_1,\,i_2\,)}g_N^{\emptyset}\biggr|_1\\
\!\!\!&\le \biggl|\sum_{m=N+1}^{\infty}\lambda_m
\,\biggl(\frac{1}{\sqrt{n}} \sum_{0 \le\, i \, \le n-1}\varphi_m
\circ T^{i}\biggr)\otimes \biggl( \frac{1}{\sqrt{n}}\sum_{0 \le \,
i \, \le n-1} \varphi_m \circ T^{i} \biggr)\biggr|_{\,2,\,2,\,\pi}\\
&\le \sum_{m=N+1}^{\infty}|\lambda_m| \underset{N\to\infty}{\to}0
\end{split}
\end{equation*}
holds uniformly in $n$ (we used here that the functions
$(\varphi_m \circ T^{i})_{1 \le m,\, 1 \le i}$ are orthonormal),
the proof is completed.
\end{proof}

\section{Exemplary applications} \label{applic}
In this section we show how the results of the present paper can
be applied in situations familiar to specialists in limit
theorems for dynamical systems or weakly dependent random
variables. We develop only a few  of all possible applications and
we do not optimize our assumptions. Instead, we show how certain
earlier known and some new results can be deduced from ours.
Applications of Theorem \ref{slln} were given in Corollaries
\ref{particular} and \ref{sllnser}.

\subsection{Doubling transformation}
 Let $X=\{z \in \C: |z|=1\},$ $\mu$ be the probability Haar
measure on $X,$ $Tz=z^2, z \in X.$  Clearly,
$$(Vf) (x) = f (x^2), \,\,\,  (V^*f)(x)=1/2\sum_{\{u:\,u^2=x\}} f(u).$$
\vspace{-0.1cm} $T$ is known to be exact \cite{Ro1961}.
 If $ f_1 \in L^2(\mu)$ and $ \int_{X}f_1(x) \mu(dx)=0$
then the series
$$ \sum_{k \ge 0} V^{*\,k} f_1 $$
converges in $L^2(\mu) $ under very mild conditions. For example,
the condition
$$ \sum_{k \ge 0}w^{(2)}(f_1,2^{-k}) < \infty$$
is sufficient. Here $w^{(2)}(f_1,\cdot)$ is the modulus of
continuity of $f_1$ in $L^2(\mu).$ \\

a) {\bf  Translation-invariant kernels.} Let now $f \in
L^2(\mu^2)$ be of the form
\begin{equation} \label{diff}
f(x_1,x_2) = g(x_1 x_2^{-1})
\end{equation}
 with some
$ g(x)= \sum_{k \in \Z} g_k x^k \in L^2(\mu) $. Assume that $f =
f_2$ (that is $f$ is canonical), real-valued and symmetric. This
means that $g_0=0,$ $ g_k$ are real and satisfy $g_{-k}=g_k$ for
all $ k \in \Z.$ Assume, moreover, that $f_2 \in
L_{2,\,\pi}^{sym}(\mu^2).$ In our setup this is equivalent  to the
relation
\begin{equation} \label{absconv}
\sum_{k \in \Z } |g_k| < \infty.
\end{equation}
The condition of the existence of the limit
$$\lim_{n \to \infty}\sum_{0\, \le \, i_1, \, i_2 \, \le \, n-1 }
V^{* (\,i_1,\,i_2)}f_2 $$ in $ L_{2,\,\pi}(\,\mu^2)$ is satisfied
if the series $\sum_{k=0}^{\infty} nV^{*n}g $ is norm convergent
in the space of absolutely convergent trigonometric series, that
is
$$ \sum_{k \in \Z} \sum_{n \ge 0} n |\,g_{2^n k}| < \infty .$$
The latter condition holds, for example, if for some $C>0$  and
$\d>0$
$$|g_m| \le \frac{C}{|m |(\log |m|)^{1+\d}}$$
for every $ m \in \Z, m \neq 0$. This condition is a very mild
strengthening of \eqref{absconv}; in its turn, the latter is, for
kernels $f$ of the type \eqref{diff}, a necessary and sufficient
condition
to belong to $L_{2,\,\pi}(\,\mu^2)$.\\
 b) {\bf General kernels.} Consider now
(compare Proposition \ref{suffic}) a general kernel $f$ $\in$
$L_2(X^2,\mathcal F^{\otimes 2}, \mu^2)$ with Fourier expansion
$$f (x_1,x_2)=\sum_{k_1,\,k_2 \in \,\Z} f_{k_1,\,k_2}x_1^{k_1}
x_2^{k_2}, \,\,\,\,\,x_1,x_2 \in X.$$ Assume that the kernel  $f$
is real-valued and  symmetric, that is $f_{-k_1,\,-k_2}=\overline
f_{k_1,\,k_2}$ and $f_{k_2,\,k_1}=f_{k_1,\,k_2}$ for $k_1,\,k_2
\in \Z$. Following the notation of Remark \ref{4.3}, we have
$f_0=f_{0,\,0}$, $f_1(x)=\sum_{k \in \Z\setminus \{0\}} f_{k,\,0}
\,x^k$, $f_2(x_1,x_2)=\sum_{k_1,\,k_2 \in\, \Z \setminus \{0\}}
f_{k_1,\,k_2}\, x_1^{k_1} x_2^{k_2}$. The kernel $f$ satisfies all
conditions of Theorems 2 and 4 whenever
$$\sum_{n \ge0}\Bigl(\sum_{k \in  \,\Z\setminus \{0\}} |f_{2^nk,\,0}|^2_2\Bigr)^{1/2}< \infty
\,\,\, \text{and} \,\,\, \sum_{n_1,\,n_2 \ge \, 0}\,\,
\sum_{k_1,\,k_2 \in  \, \,\Z\setminus \{0\}}
|f_{2^{n_1}k_1,\,2^{n_2}k_2}|\,< \infty.$$
\begin{remark} \label{expand}
In this subsection we gave applications of our results to the
simplest example of a \emph{differentiable expanding map}.
This is based on the group structure
of the example and its Fourier
analysis.
A more general approach can be developed on the basis of
the \emph{transfer operator} ($V^*$ in our setting)
restricted to some spaces of nice (smooth, H\"older or Sobolev)
functions.
\end{remark}
\subsection{Stationary processes (martingale kernels, mixing conditions, Markov
processes)} \label{stat}
 Let $\xi=(\xi_n)_{n \in \Z}$ be an ergodic stationary random process
defined on the space $(X, \mathcal F, \mu)$ where an invertible
measure preserving transformation $T$ acts so that
$\xi_{n+1}=\xi_n\circ T,\, n \in \Z$. We assume that all $\xi_n$
take values in a probability space $ (Y, \mathcal G, \nu)$, $\nu$
being the common distribution of $(\xi_n)_{n \in\, \Z}$.  Let
$(X^d, \mathcal F^{\otimes d}, \mu^d)$ be the $d$-th Cartesian
power of $(X, \mathcal F, \mu)$ with the coordinatewise action of
$(T^{\mathbf n})_{\mathbf n \in\, \Z^d}$ and the corresponding
operators $(V^{\mathbf n})_{\mathbf n \in\, \Z^d}$; let,
furthermore, $(\xi_n^{(i)})_{n \in \Z}$, $1 \le i \le d$, be
independent copies of $(\xi_n)_{n \in \Z}$ defined on $(X^d,
\mathcal F^{\otimes d}, \mu^d)$ so that $\xi_n^{(i)}
(x_1,\dots,x_d)=\xi_n(x_i)$, where $\,x_1,\dots,$ $ x_d \in X$,
$1\le i \le d, n \in \Z.$ Assume now that we are given some $F \in
L_{p,\,\pi}(Y^d,\mathcal G^{\otimes d}, \nu^d) $ for some $d \in
\N$ and $p \in [1,\infty)$. Then $f=F(\xi_0^{(1)},
\dots,\xi_0^{(d)}) \in L_{p,\,\pi}(X^d,\mathcal F^{\otimes d},
\mu^d)$, $F(\xi_{n_1}^{(1)},\dots,\xi_{n_d}^{(d)})=V^{\mathbf n}f$
and\\ $F(\xi_{n_1},\dots,\xi_{n_d})=D_d V^{\mathbf n}f$  for every
$\mathbf n=(n_1, \dots, n_d) \in
\Z^d$.\\
In the rest of the paper, instead of saying that  an assertion of the previous part of the paper applies to a kernel $f$ and a  transformation $T$, we will usually say that  this assertion applies to the kernel $F$ (the process $\xi$ will be omitted).

\ref{stat}.1. {\bf Martingale kernels.} Let $d=2$. Set $\mathcal
F_0=\s(\xi_0,\xi_{-1},\dots),$ the $\s$-field generated by
$\xi_0,\xi_{-1},\dots$, and $\mathcal F_n = T^{-n}\mathcal
F_0=\s(\xi_n,\xi_{n-1},\dots)$. Assume that $f=
F(\xi_0^{(1)},\xi_0^{(2)})$ is a canonical kernel. Obviously, it
is measurable with respect to $\mathcal F_0^{(1)} \otimes \mathcal
F_0^{(2)}\bigl(\,\overset{\mathrm{def}}{=}
\s((\xi_0^{(1)},\xi_{-1}^{(1)},\dots,\xi_0^{(2)},\xi_{-1}^{(2)},\dots)\bigr)$.

 The equivalent of
\eqref{conv8} for
invertible $T$ is the existence
of the limit
$$
\underset{n_1,\,n_2 \to \infty}{\lim}\underset{\substack{0 \, \le\, i_1 \,\le \, n_1-1\\
0\, \le \,i_2\, \le \, \, n_2-1}}{\sum}
V^{(i_1,\,i_2)}E\bigl(F(\xi_0^{(1)},\xi_0^{(2)})|\,\mathcal
F_{-i_1}^{(1)} \otimes \mathcal F_{-i_2}^{(2)}\bigr)
$$
in the space $L_{2,\,\pi}(X^2,\mathcal F^{\otimes 2}, \mu^2)$. To
compare our Theorem \ref{degenCLT} (in its invertible
modification, see Remark \ref{adapt})  with the main limit theorem
in \cite{LeNe2011} notice that it is assumed there that the kernel
$F$ is \emph{symmetric} and satisfies $
E\bigl(F(\xi_0^{(1)}\!,\,\xi_0^{(2)})|\,\mathcal F_{-1}^{(1)} \otimes
\mathcal F_{0}^{(2)}\bigr)=0. $ This implies that a non-vanishing
summand may appear in the above sum only for
$i_1=i_2=0$, so we have nothing more to check in this case.\\

\ref{stat}.2. {\bf Processes satisfying mixing conditions.}  For
$k\in \Z$ we set $\mathcal F_k=\s(\xi_l, l \le k)$, $\mathcal
F^k=\s(\xi_l, l \ge k)$, $\mathcal F(k)=\s(\xi_k)$;
let $E_k,E^k,E(k)$  denote the corresponding conditional
expectation operators, $E$ being the unconditional expectation.
For the system of $\s$-fields $(\mathcal F_k, \mathcal
F^k)_{k \in\, \Z}$ and $n \in \Z_+$ define the well-known \emph{mixing
coefficients} by setting
$$ \a(n)=\underset{A \in\, \mathcal F_k, \, B \in \,\mathcal F^{k+n}}{\sup}|\,\mu(A \cap B) - \mu(A)
\mu(B)|,$$
$$ \varphi(n)=\underset{A \in \,\mathcal F_k,\,\,\mu(A)
>\,0, \, B \in\, \mathcal F^{k+n}}{\sup}\mu^{-1}(A)|\,\mu(A \cap B) -
\mu(A) \mu(B)|,$$
$$ \psi(n)=\underset{A \in \,\!\mathcal F_k\!,\, B \in\,\, \mathcal
F^{k+n} ,\,\,\, \mu(A)\,\mu(B)
>\,0}{\sup}\mu^{-1}(A)\mu^{-1}(B)|\,\mu(A \cap B) -
\mu(A) \mu(B)|.$$ For the norms  of the operators $E^{k+n}E_k-E,
E_k E^{k+n}-E$ which act from $L_q(\,\mu)$ to $L_p(\,\mu)$ $(p,\,q
\in [1,\,\infty]\,) $ certain bounds  in terms of the mixing
coefficients are known \cite{Da1968,IbLi1971}. Indeed, we have for
$1\, \le p\,\le q\le\, \infty $
\begin{equation}\label{opalpha}
 \max\,(|E^{k+n}E_k-E|_{\,q,\,p}\,,|E_k E^{k+n}-E|_{\,q,\,p})\le
C(q,\,p)\, \a(n)^{\,p^{-1}-q^{-1}} \,\,\,
\end{equation}
\begin{equation}\label{opphi2}
|E_k E^{k+n}-E|_{\,q,\, p} \le 2 \varphi(n)^{1-q^{-1}},
\end{equation}
and for  $1\,\le p\,,q  \le \infty $
\begin{equation}\label{opsi}
\max\,(|E^{k+n}E_k-E|_{\,q,\,p}\,,\,|E_k E^{k+n}-E|_{q,\,p}) \le
\psi(n).
\end{equation}
Notice that if at least one of the mixing coefficients tends to 0,
the process $\xi$ is Kolmogorov (see Remark \ref{adapt}) and 
consequently ergodic. Set
$$M_{q,\,p}=\sum_{n \ge\, 0}|E_k E^{k+n}-E|_{\,q,\, p}\,,\,\, M_{q,\,p}'=\sum_{n \ge\, 0}|E^{k+n}E_k-E|_{\,q,\, p}$$
(in view of stationarity $M_{q,\,p}$ and $M_{q,\,p}'$ do not
depend on $k$).\\
In the rest of \ref{stat}.2 we show how Proposition \ref{suffic}
(more precisely, its analogue for an invertible $T$) can be
used applying the results of the paper to $V$-statistics of a
process $\xi$ with suitable mixing
properties.\\
 Let $(\e_k)_{k=0}^{\infty}$ be a sequence
of functions satisfying
\begin{equation} \label{base}
\begin{split}
&\e_k \in L_q\,(Y,\mathcal G, \nu), \,\,|\e_k|_{\,q}=1\,\,\, (k
\ge
0),\\
&\e_0\equiv 1, \int_Y \e_k(y) \nu(dx)=0 \,\,\, (k \ge 1).
\end{split}
\end{equation}
 Set $e_k=\e_k\circ\xi_0, k \ge 0,$ and fix some $p \in
[1,\,q\,].$ Observe that for every $k \ge 1$
\begin{equation} \label{basconv}
C_{p,\,k}=\!\sum_{n \ge 0} |\,E_{-n} e_k|_{\,p}=\!\sum_{n \ge 0}
|\,(E_{-n} E^0 -E)e_k|_{\,p} \le M_{q,\,p}
|\,\e_k|_{\,q}=M_{q,\,p}\,.
\end{equation}
Assume that a function $F \in L_q\,(Y^m,\mathcal G^{\otimes m}, \nu^m)$ expands
   into the series
\begin{equation} \label{multiple1}
F(y_1,\dots,y_m)=
\sum_{\mathbf 0 \boldsymbol{<} \mathbf k
\boldsymbol{<}\boldsymbol{\infty}}
\l^F_{\,\bold k}\, \e_{k_1}(y_1)\cdots \e_{k_m}(y_m)
\end{equation}
with a
some family $(\l^F_{\mathbf k})_{\,\mathbf 0
\boldsymbol{<}\mathbf k \boldsymbol{<}\boldsymbol{\infty}}$. For
the following  expression of the type of \eqref{up} we have
\begin{equation} \label{up1}
 C_p^F\eqdef\sum_{\mathbf 0
\boldsymbol{<}\mathbf k \boldsymbol{<}\boldsymbol{\infty}}
|\,\l^F_{\,\mathbf k}\,|\, C_{p,\,k_1}\cdots\, C_{p,\,k_m}\le
(M_{q,\,p})^m \sum_{\mathbf 0 \boldsymbol{<}\mathbf k
\boldsymbol{<}\boldsymbol{\infty}} |\,\l^F_{\,\mathbf k}\,|\,,
\end{equation}
so $ C_p^F < \infty$  whenever
$$M_{q,\,p} <\infty$$
(this is a condition on the mixing rate of the process $\xi$) and
expansion \eqref{multiple1} of the function $F$ satisfies the
condition
\begin{equation} \label{coeff} \sum_{\mathbf 0
\boldsymbol{<}\mathbf k \boldsymbol{<}\boldsymbol{\infty}}
|\,\l^F_{\,\mathbf k}\,|\, < \infty.
\end{equation}
Thus,\emph{ the invertible version of Proposition \ref{suffic}
applies to the kernel $f:(x_1,\dots,x_m) \mapsto
F(\xi_0(x_1),\dots,\xi_0(x_m))$ with some $p \in [1,\infty]$ and
the system $(e_k)_{k=0}^{\infty}$ if, for a certain $q \in
[p,\infty]$, the system $(\e_k)_{k=0}^{\infty}$ satisfies the
conditions \eqref{base},  $F \in L_{\,q}(\,Y^m,\mathcal G^{\otimes
m}, \nu^m) $ admits the representation \eqref{multiple1}, 
satisfying \eqref{coeff}, and we have
$M_{q,\,p}< \infty$ for the process $\xi$}.\\
 We now indicate
conditions (stated in terms of $\alpha, \varphi$ and $\psi$)
 under which Theorems 2, 3 and 4 of the paper,
in their invertible forms
 and numerated by 2$^{\,\prime}$, 3$^{\,\prime}$ and 4$^{\,\prime}$,  apply to an $F$. Theorem 3
 needs more substantial changes in case of the mixing coefficient $\varphi$.
Below $(\e_k)_{k \ge 0}$
is a system satisfying \eqref{base} with some parameter $q$.\\

{\bf a)} Let  $q\in [2d,\infty]$. We will use
\eqref{opalpha}, \eqref{opphi2} and \eqref{opsi}, substituting
there, in place of the pair $(q,p)$, the pair 
$(q,2d)$; we will employ Proposition \ref{suffic} and formulas \eqref{basconv}, \eqref{up1} with $p=2d$. Theorem \ref{clt2}$^{\,\prime}$ applies to
an $F \in L_{2}^{sym}(\,\nu^d) $ if
\begin{trivlist}
\item 1)\,\,\,\,\, at least one of the series
\begin{equation} \label{ser1}
\sum_{n \ge 0} \alpha(n)^{\!
(2d)^{-1}\!-q^{-1}},\, \sum_{n \ge
0}\varphi(n)^{1-\!
q^{-1}\!},\, \sum_{n \ge 0} \psi(n) \end{equation}
converges (for $d=q$ the convergence of the $\alpha$-series means
that $\alpha(n)=0$ for $n \ge n_0$), and
 \item 2) \,\,\, for every $m =2,
\dots,d$ $R_m F$ belongs to $L_q^{sym}(\nu ^m)$ and admits the
representation
\begin{equation} \label{decomp}
R_m F(y_1,\dots,y_m)=
\sum_{\mathbf 0 \boldsymbol{<} \mathbf k
\boldsymbol{<}\boldsymbol{\infty}}
\l^{R_m F}_{\,\bold k}\, \e_{k_1}(y_1)\cdots \e_{k_m}(y_m)
\end{equation}
where the coefficients satisfy $ \sum_{\,\mathbf 0 \,
\boldsymbol{<}\mathbf k \, \boldsymbol{<}\,\boldsymbol{\infty}}
|\,\l^{R_m F}_{\,\mathbf k}\,|\, < \infty $.\\
Under condition 2) with $q = 2d$ Theorem \ref{clt2}$^{\,\prime}$
applies, in particular, if  $ \sum_{n \ge
0}\varphi(n)^{1-\!
(2d)^{-1}\!} < \infty$.
\end{trivlist}
\vspace{0.23cm}

 {\bf b)}\,\,\, To simplify the
 statements involving $\varphi$ assume
that $d \ge 2$. Let $q\in [d,\infty]$.  
We will use
\eqref{opalpha}, \eqref{opphi2} and \eqref{opsi}, substituting
there, in place of the pair $(q,p)$, the pair 
$(q,d)$; we will employ Proposition \ref{suffic} and formulas \eqref{basconv}, \eqref{up1} with $p=d$. 
 Theorem \ref{clt1}$^{\,\prime}$ applies to
an $F \in L_{1}^{sym}(\,\nu^d) $ if
\begin{trivlist}
\item 1)\,\,\,\,\, at least one of the series
\begin{equation} \label{ser1}
\sum_{n \ge 0} \alpha(n)^{d^{-1}\! -\! q^{-1}}\, \sum_{n \ge
0}\varphi(n)^{1-q^{-1}},\, \sum_{n \ge 0} \psi(n) \end{equation}
converges (if $q=d$ the convergence of the $\alpha$-series means
that $\alpha(n)=0$ for $n \ge n_0$); \item 2) $R_1 F$ satisfies
the relation \eqref{sum_1}:
$$ \Bigl|\sum_{k=0}^{n-1}(R_1 F)\circ\xi_k\Bigr|_1\,\,=\,\, O(\sqrt{n})\,\,;$$
 \item 3) \,\,\, for every $m =2,
\dots,d$ $R_m F$ belongs to $L_q^{sym}(\nu ^m)$ and admits the
representation
\begin{equation} \label{decomp}
R_m F(y_1,\dots,y_m)=
\sum_{\mathbf 0 \boldsymbol{<} \mathbf k
\boldsymbol{<}\boldsymbol{\infty}}
\l^{R_m F}_{\,\bold k}\, \e_{k_1}(y_1)\cdots \e_{k_m}(y_m)
\end{equation}
where the coefficients satisfy $ \sum_{\,\mathbf 0 \,
\boldsymbol{<}\mathbf k \, \boldsymbol{<}\,\boldsymbol{\infty}}
|\,\l^{R_m F}_{\,\mathbf k}\,|\, < \infty $.\\
Under conditions 2) and 3) Theorem \ref{clt1}$^{\,\prime}$
applies, in particular, if $ q=2d$ and $ \sum_{n\, \ge \,
0}\alpha(n)^{1/2d} < \infty $.
\end{trivlist}
 {\bf c)} Theorem \ref{degenCLT} leads to a result on mixing
processes in the following way. Let $F \in
L^{sym}_{2,\,\pi}(Y^2,\mathcal G^{\otimes 2}, \nu^2) $ be a
canonical function. Hence, it is the kernel of a \emph{nuclear}
(or \emph{trace class}) symmetric integral operator in $L_2(\nu)$
vanishing on constant functions. The general theory says that in
$L_2(\nu)$ there exists an orthogonal normalized sequence
$\e_0\equiv 1,\e_1,\dots$ and a real sequence $\g_1, \g_2, \dots$
such that
\begin{equation} \label{eigen}
 F(x_1,x_2)=\sum_{k=1}^{\infty} \g_k \,
\e_k(x_1)\,\e_k(x_2),
\end{equation}
 where $\sum_{k=1}^{\infty}|\,\g_k| < \infty$ ($k=0$ is omitted
 because $F$ is canonical). Neglecting
the assumption of canonicity and symmetry, such functions form
exactly the space $L_{2, \, \pi}(\nu^2)$; the projective norm
agrees for symmetric functions with the sum of moduli of the
eigenvalues of the corresponding integral operators. Thus
$f:(x_1,x_2) \mapsto F(\xi_0(x_1),\xi_0(x_2))$ is a function to
apply Proposition \ref{suffic} with $d=2$ and $e_k=\e_k \circ\xi_0
\,(k \ge 0)$. Then for $k \ge 1$ $C_{2,\,k} \le M_{2,\,2}$. The
latter quantity is bounded above by any of the series $\sum_{n \ge
0}\,\varphi(n)^{1/2}$, $\sum_{n \ge 0}\, \psi(n)$. Thus, the
invertible version of Theorem \ref{degenCLT} applies whenever at
least one of these series converges.
\begin{remark} \label{asymm} The last assertion under the
assumption $\sum_{n \ge 0}\,\varphi(n)^{1/2} < \infty $ is, up to
inessential details, Theorem 5 in \cite{Ea1979}. In
\cite{BoVo2008} the authors express their doubts on correctness in
\cite{Ea1979} to substituting a dependent process into the
function \eqref{eigen}. Our conclusion agrees with that of
\cite{Ea1979}. In our paper the correctness is a simple
consequence of general properties of projective tensor products.
However, an elementary reasoning shows that the series
\eqref{eigen}  absolutely converges in $L_1(X^2,\kappa)$ where
$\kappa$ is an arbitrary probability on $X^2$ with one-dimensional
marginals $\mu$.
\end{remark}
\ref{stat}.3. {\bf Discrete time Markov processes.}
 Let $\xi=(\xi_n)_{n \, \in \, \Z}$ be a stationary Markov process
defined on the space $(X, \mathcal F, \mu)$ where an invertible
measure preserving transformation $T$ acts so that
$\xi_{n+1}=\xi_n\circ T,\, n \in \Z$.
We assume that all $\xi_n$ take values in a probability space $
(Y, \mathcal G, \nu)$, $Y$ being the \emph{ state space} of $\xi$
and $\nu$ its \emph{stationary distribution}. We will use
the notations $\mathcal F_k$, $\mathcal F^k$, $\mathcal
F(k)$, $E_k,E^k, E(k)$ and $E$ as introduced above.\\
 Let $Q$ be the \emph{transition
operator} of $\xi$ acting on every space $L_p(\nu),$ $1 \le p \le
\infty,$ with norm $1$ and satisfying $E_k
f(\xi_{k+1})=(Qf)(\xi_k)$ for every $f \in L_1(\nu)$ and $k \in
\Z$. Assuming $\mathcal F =\s(\xi_l, l \in \Z)$, the process $\xi$
(that is the transformation $T$) is ergodic if and only if for the
transition operator $Q: L_2(\nu) \to L_2(\nu)$ every solution to
the equation $Qf=f$ is a constant. To stay within the assumptions
of the present paper we assume a stronger relation $Q^n h
\underset{n \to
\infty}{\to} \int h (y)\nu(dy)$ $(h \in L_1(\nu))$ which implies the Kolmogorov property of $\xi$.\\
Let $d \ge 1$ and $(\e_k)_{k=0}^{\infty}$ be a sequence of
functions satisfying \eqref{base} with $q=2d$. Let $I_{\nu}$
denote the identity operator in every space $L_{q}(\nu)$. Assume that for
some $C > 0$ and every $k \ge 1$ the equation
$(I_{\nu}-Q)\phi_k=\e_k$ is solvable and $|\,\phi_k\,|_{\,2d}\le C
$ (notice that the latter condition is fulfilled if the
restriction $(I_{\nu}-Q)|_{L^0_{2d}}$ is invertible, $L^0_{2d}$
denoting the subspace of functions in $L_{2d}$ with integral $0$).
Let  $F\in L_{2}\,(Y^d,\mathcal G^{\otimes d}, \nu^d)$
satisfy assumption 2) of paragraph {\bf a)} in 9.2.2 with $q=2d$. Let, finally, 
the equation $(I_{\nu}-Q)g=R_1F$ have a solution $g \in L_2(\nu)$.  Then 
 Theorem \ref{clt2}$^{\,\prime}$ applies to $f=F(\xi_0^{(1)}, \dots,\xi_0^{(d)}) $.\\

\section{Acknowledgment}
The authors would like to thank Herold Dehling for several
discussions clarifying  many aspects of limit distributions for
$V$-statistics and for his encouragement to get this paper
written. Also comments by referees were very helpful.
 The research was supported by the Deutsche
Forschungsgemeinschaft under Grant Number 436 RUS 113/962/0-1 and
the Russian Foundation
for Basic Research under Grant Number 09-01-91331-NNIO-a.\\
Manfred Denker was also partially supported by the National
Science Foundation
under Grant Number DMS-10008538.\\
Mikhail Gordin was also partially supported by the Grant Number
13-01-00256-a of the Russian Foundation for Basic Research   and
by the grant Number NS-2504.2014.1   for the Support of Scientific Schools. He thanks Axel Munk (Institute for Mathematical
Stochastics) and Laurent Bartholdi (Mathematical Institute) for
their hospitality at the University of G\"ottingen where a part of
this paper was prepared.


\begin{thebibliography}{3}

\bibitem{Aa96}
 Aaronson J., Burton R., Dehling H., Gilat D., Hill T., Weiss B.
 \emph{Strong laws for $L$- and $U$-statistics.} Trans. Amer. Math. Soc. {\bf 348}: 7 (1996),
  2845--2866.

\bibitem{ArGi1993}
Arcones M. A., Gin{\'e} E. \emph{Limit theorems for
$U$-processes}. Ann. Probab. {\bf 21}(1993): 3, 1494--1542.


\bibitem{Ama} Amanov A. K.
\emph{Limit distribution of a von Mises functional with a
degenerate kernel for dependent random variables.} Probabilistic
models and mathematical statistics (Russian), 7--14, no.173,
``Fan'', Tashkent, 1987.

\bibitem{Bab1} Babbel B. \emph{
Schwache Invarianzprinzipien f\"ur verallgemeinerte von Mises
Funktionale und U-Statistiken von Funktionalen schwach abh\"angiger
Prozesse im Mehrstichprobenfall.} (Weak invariance principles for
generalized von Mises functionals and U-statistics of functionals
of weakly dependent processes in the multisample case) (German).
G\"ottingen (FRG), 1987, Univ. G\"ottingen, Mathematisch-Naturwissenschaftlicher Fachbereich, Diss.

\bibitem{Bab1989} Babbel B. \emph{Invariance principles for $U$-statistics and von Mises
functionals.}  J. Statist. Plann. Inference  {\bf 22}: 3 (1989),
337--354.

\bibitem{Bi1961} Billingsley P. \emph{The Lindeberg--Levy theorem for martingales}. Proc.
Amer. Math. Soc.\,{\bf 12} (1961), 788--792.

\bibitem{BoBy2006-1} Borisov I. S., Bystrov A. A. \emph{Stochastic integrals
and asymptotic analysis of canonical von Mises statistics based on
dependent observations.}  High dimensional probability,  1--17,
IMS Lecture Notes Monogr. Ser., {\bf 51}, Inst. Math. Statist.,
Beachwood, OH, 2006.

\bibitem{BoBy2006-2} Borisov I. S., Bystrov A. A. \emph{
Limit theorems for canonical von Mises statistics constructed from
dependent observations.} (Russian). Sibirsk. Mat. Zh.\,{\bf 47}:
6, (2006), 1205--1217; transl. in Siberian Math. J.\,{\bf 47}: 6,
(2006), 980--989.

\bibitem{BoVo2008} Borisov I. S., Volod'ko N. V.
\emph{Orthogonal series and limit theorems for canonical $U$- and
$V$-statistics of stationarily connected observations.} (Russian)
Mat. Tr.\,{\bf 11}: 1 (2008), 25--48; transl. in Siberian Adv. in
Math. {\bf 18}: 4 (2008), 242--257.


\bibitem{BoBuDe2001} Borovkova S., Burton R., Dehling H.
\emph{Limit theorems for functionals of mixing processes with
applications to $U$-statistics and dimension estimation}. Trans.
Amer. Math. Soc. {\bf 353}: 11 (2001), 4261--4318 (electronic)

\bibitem{BoBuDe2002} Borovkova S., Burton R., Dehling H.
 \emph{From dimension estimation to asymptotics of dependent
$U$-statistics}. In \emph{Limit theorems in probability and
statistics, Vol. I (Balatonlelle, 1999)}.   J\'anos Bolyai Math.
Soc., Budapest, 2002, 201--234.

\bibitem{Brad1988}Bradley R. C.
\emph{On some results of M. I. Gordin: a clarification of a
   misunderstanding}. J. Theoret. Probab. {\bf 1}: 2 (1988), 115--119.

\bibitem{Burk1966}Burkholder D. L.
\emph{Martingale transforms}. Ann. Math. Statist. {\bf 37} (1966),
1494--1504.

\bibitem{Da1968} Davydov Yu. A.  \emph{Convergence of distributions generated
by stationary stochastic processes.} (Russian). Prob. Theory and
its Appl. {\bf 13}: 4 (1968), 730--737.



\bibitem{DeDePhi1984}
Dehling H., Denker M., Philipp W. \emph{Invariance principles for
von Mises and $U$--statistics.} Z. Wahrsch. Verw. Gebiete {\bf 67}: 2 (1984), 139--167.


\bibitem{DeFl1993} Defant, A., Floret, K. \emph{Tensor norms and
operator ideals.} North-Holland Mathematics Studies {\bf 176},
North-Holland, Amsterdam -- London -- New York -- Tokyo, 1993.


\bibitem{DeGi1999} de la Pe{\~n}a V. H., Gin{\'e}
E. \emph{Decoupling. From dependence to independence.
   Randomly stopped processes. $U$--statistics and processes. Martingales and
   beyond}. Probability and its Applications (New York),
   Springer-Verlag, New York, 1999.


\bibitem{DeTa} Dehling H., Taqqu M. S. \emph{The limit behavior of empirical processes
and symmetric statistics for stationary sequences.} Proceedings of
the 46th Session of the International Statistical Institute, Vol.
4 (Tokyo, 1987). Bull. Inst. Internat. Statist. {\bf 52}: 4
(1987), 217--234.

\bibitem{De} Dehling H. \emph{Limit theorems for dependent $U$-statistics.}  Dependence
in probability and statistics,  65--86, Lecture Notes in Statist., 187, Springer,
New York, 2006.

\bibitem{DeWe} Dehling H., Wendler M. \emph{Central limit theorem and a
    bootstrap for $U$-statistics of strongly mixing data.} J. Mult. Anal. {\bf
    101} (2010), 126--137.

\bibitem{DeKe84} Denker M., Keller G. \emph{On $U$-statistics and v. Mises' statistics
for weakly dependent processes.}  Z. Wahrsch. Verw. Gebiete  {\bf
64}: 4 (1983),  505--522.

\bibitem{DeKe86} Denker M., Keller G. {\em Rigorous statistical procedures for data from
dynamical systems.}  J. Statist. Phys.  {\bf 44}: 1-2 (1986),
67--93.

\bibitem{DeGrKe85}  Denker M., Grillenberger C., Keller, G. {\em A note on invariance principles
for v. Mises' statistics.}  Metrika  {\bf 32}: 3-4 (1985),
197--214.

\bibitem{DuSchw1958}
Dunford N., Schwartz J. T., \emph{Linear Operators. I. General
Theory}, Interscience Publishers, Inc., N. Y., 1958.


\bibitem{DyMa1983} Dynkin E.B., Mandelbaum A. \emph{Symmetric
statistics, Poisson point proceses, and multiple Wiener
integrals}. Ann. of Stat. {\bf 11}: 3 (1983), 739--746.




\bibitem{Ea1979} Eagleson G. K.\emph{ Orthogonal expansions and
$U$-statistics}. Austral. J. Statist.{\bf 21}: 3 (1979), 221--237.


\bibitem{Fil1962}  Filippova A. A. \emph{Mises' theorem on the asymptotic behavior
    of functionals of empirical distribution functions and its statistical
    applications.} Theory of Probability and Applications {\bf 7} (1962), 24--57.



\bibitem{GoKr1969} Gohberg I. C.,  Krein M. G.  \emph{Introduction to the theory of
linear nonselfadjoint operators}. Translated from the Russian 1965
original. Translations of
 Mathematical Monographs, Vol. {\bf 18} American Mathematical Society, Providence,
 R.I., 1969.



\bibitem{GoKrGo2000} Gohberg I., Krupnik N., Goldberg S.
\emph{Traces and determinants of linear operators.} Operator
Theory, Advances and Applications {\bf 116}. Birkh\"auser, 2000.


\bibitem{Go1969}  Gordin M. I. \emph{On the central limit theorem for
stationary processes.}  (Russian)   Dokl. Akad. Nauk SSSR,
\textbf{188} (1969), 739--741. Transl.: Soviet Math. Dokl.
\textbf{10} (1969), 1174--1176.

\bibitem{Go1973}  Gordin M. I. \emph{Central limit theorem for
stationary processes without the assumption of finite variance}.
(Russian). Abstracts of Communications, T. 1: A-K, p.p. 173--174.
International Conference on Probability Theory and Mathematical
Statistics, June 25--30, 1973, Vilnius.

\bibitem{GoLi1978} Gordin M. I., Lifshits B. A.  \emph{Central limit
theorem for stationary Markov processes.} (Russian) Dokl. Akad.
Nauk SSSR  \textbf{239} (1978), 766--767; transl.: Soviet Math.
Dokl. \textbf{19}: 2 (1978), 392--394.

\bibitem{Go2009} Gordin M. I. \emph{Martingale-coboundary representation for a class of random
fields}. (Russian). Zap. Nauchn. Sem. S.-Peterburg. Otdel. Mat.
Inst. Steklov (POMI), \textbf{364} (2009), Veroyatnost i
Statistika. no. 14.2, pp.88--108, 236; transl.: Journ. Math. Sci.
(NY) {\bf 163}: 2 (2009), 363--374.

\bibitem{GraPro} Grassberger P., Procaccia I.
\emph{Characterization of
    strange attractors.}  Phys. Reviews Letters  {\bf 50}  (1983),  346--349.

\bibitem{Hal1946} Halmos P. \emph{The theory of unbiased estimation.}
  Ann. Math. Statistics {\bf 17} (1946), 34--43.

\bibitem{Hoe1948} Hoeffding W.
\emph{A class of statistics with asymptotically normal distribution.}
Ann. Math. Statistics  {\bf 19}  (1948). 293--325.


\bibitem{Ib1963} Ibragimov I. A. \emph{A central limit theorem for a
class of dependent random variables}. (Russian), Prob. Theory and
its Applic. {\bf 8}: 1 (1963), 89--94.

\bibitem{IbLi1971} Ibragimov I. A., Linnik Yu. V.
\emph{Independent and stationry sequences of random variables.}
Wolters-Noordhof Publishing, Groningen, 1971.



\bibitem{Kha} Khashimov Sh. A.
\emph{Asymptotic normality of a generalized von Mises functional
for $m$-dependent variables.} Theory of random processes and its applications
(Russian), 141--148, ``Naukova Dumka'', Kiev, 1990.

\bibitem{KoBo1994} Koroljuk V. S., Borovskich Yu. V. \emph{Theory of
$U$-statistics.} Mathematics and its Applications, {\bf 273}.
Translated from the 1989 Russian original. Kluwer Ac. Publ. Group,
Dordrecht.

\bibitem{LeNe2011} Leucht A., Neumann M. H. \emph{Degenerate U- and V-statistics under
ergodicity: asymptotics, bootstrap and applications in
statistics}. Ann. Inst. Statist.Math. \textbf{65} (2013), 349--386.

\bibitem{LeNe2012} Leucht A., Neumann M. H. \emph{Dependent wild
bootstrap for degenerate U- and V-statistics}. J. Multivariate Analysis \textbf{117}  (2013), 257--280.

\bibitem{Le2012} Leucht A. \emph{Degenerate U- and V-statistics under weak
dependence: Asymptotic theory and bootstrap consistency}.
Bernoulli \textbf{18}: 4 (2012), 552--585.

\bibitem{Ma1978} Maigret N. \emph{ Th\'eor\`eme de limite centrale fonctionnel pour une cha\^\i ne
de Mar\-kov r\'ecurrente au sens de Harris et positive.} Ann.
Inst. H. Poincar\'e Sect. B (N.S.), \textbf{14}: 4 (1978),
425--440.



\bibitem{Maj2005} Major P. \emph{Tail behaviour of multiple random integrals and
$U$-statistics}. J Probab. Surv. {\bf 2} (2005), 448--505.


\bibitem{Neu} Neuhaus G. \emph{Functional limit theorems for $U$-statistics in the
degenerate case.}  J. Multivariate Anal.  {\bf 7}: 3 (1977),
424--439.


\bibitem{ReSi1980} Reed M., Simon B. \emph{Methods of modern mathematical physics.
Vol.1. Functional analysis.} Academic Press, Inc., 1980.


\bibitem{Ro1961} Rohlin (Rokhlin) V. A.
\emph{Exact endomorphisms of a Lebesgue space.}
 (Russian), Izv. Akad. Nauk SSSR Ser. Mat., {\bf 25} (1961), 499--530.


\bibitem{Ry2002}Ryan R. A.
   \emph{Introduction to tensor products of Banach spaces},
   Springer-Verlag London Ltd., London, 2002.

\bibitem{Sha} Sharipov O. Sh.\emph{
The invariance principle for $U$-statistics and von Mises
functionals of weakly dependent observations.} (Russian). Teor.
Veroyatnost. i Primenen. {\bf 47} (2002), no. 4, 814--817; transl.
Theory Probab. Appl. {\bf 47}: 4 (2003), 730--733.



\bibitem{Vol2011} Volod'ko N. V.,
   \emph{Limit theorems for canonical von Mises statistics and
   $U$-statistics of $m$-dependent observations}.
   (Russian). Teor. Veroyatn. Primen.,
   {\bf 55}: 2(2010), 226--249. Transl. in Theory Probab. Appl.,\,{\bf
   55}: 2 (2011), 271--290.


\bibitem{vMis1947} von Mises R. {\em On the asymptotic distribution of differentiable statistical
functions.}  Ann. Math. Statistics  {\bf 18} (1947), 309--348.

\bibitem{Yosh} Yoshihara K.-I. \emph{Limiting behavior of $U$-statistics
for stationary, absolutely regular processes.}  Z. Wahrsch.  Verw.
Gebiete {\bf 35}: 3 (1976), 237--252.

\bibitem{Yosh92} Yoshihara K.-I. \emph{Limiting behavior of $U$-statistics
for strongly mixing sequences.}  Yokohama Math. J.  {\bf 39}: 2
(1992), 107--113.
\end{thebibliography}
\end{document}